\documentclass[11pt]{article}

\usepackage{latexsym,color,graphicx,amsthm,amsmath,amssymb}
\usepackage{lineno}

\def\reals{{\mathbb R}}
\def\cplx{{\mathbb C}}
\def\A{{\cal A}}
\def\C{{\cal C}}
\def\F{{\cal F}}

\def\XX{{1.13}}

\def\P{{\mathbb P}}
\def\O{{\mathcal O}}

\def\eps{{\varepsilon}}

\def\P{{\mathbb P}}
\newcommand{\RR}{\ensuremath{\mathbb R}}

\theoremstyle{plain}
\newtheorem{theorem}{Theorem}[section]
\newtheorem{corollary}[theorem]{Corollary}

\newtheorem{lemma}[theorem]{Lemma}

\newcommand{\Deg}{D} 
\makeatletter
\newcommand{\ProofEndBox}{{\ifhmode\unskip\nobreak\hfil\penalty50 \else
          \leavevmode\fi\quad\vadjust{}\nobreak\hfill$\Box$
            \finalhyphendemerits=0 \par}}
\makeatother

\newcommand{\proofend}{\ProofEndBox\smallskip}
\newcommand{\ignore}[1]{}
\begin{document}


\title{Incidences with curves and surfaces in three dimensions,\\ with applications to distinct and repeated distances\thanks{%
Work on this paper by Noam Solomon and Micha Sharir was supported by
Grant 892/13 from the Israel Science Foundation. Work by Micha
Sharir was also supported by Grant 2012/229 from the U.S.--Israel
Binational Science Foundation, by the Israeli Centers of Research
Excellence (I-CORE) program (Center No.~4/11), by the Blavatnik
Research Fund in Computer Science at Tel Aviv University and by the
Hermann Minkowski-MINERVA Center for Geometry at Tel Aviv
University. A preliminary version of the paper has appeared in {\it
Proc. 28th ACM-SIAM Symposium on Discrete Algorithms} (2017),
2456--2475.}}

\author{
Micha Sharir\thanks{%
School of Computer Science, Tel Aviv University, Tel Aviv 69978,
Israel. {\sl michas@tau.ac.il} } \and
Noam Solomon\thanks{%
School of Computer Science, Tel Aviv University, Tel Aviv 69978,
Israel. {\sl noam.solom@gmail.com} } }


\maketitle

\begin{abstract}
We study a wide spectrum of incidence problems involving points and curves
or points and surfaces in $\reals^3$. The current (and in fact the only viable)
approach to such problems, pioneered by Guth and Katz~\cite{GK,GK2},
requires a variety of tools from algebraic geometry, most notably
(i) the polynomial partitioning technique, and (ii) the study of algebraic surfaces
that are ruled by lines or, in more recent studies~\cite{GZ}, by algebraic curves
of some constant degree. By exploiting and refining these tools, we obtain
new and improved bounds for numerous incidence problems in $\reals^3$.

In broad terms, we consider two kinds of problems, those involving points and
constant-degree algebraic \emph{curves}, and those involving points and
constant-degree algebraic \emph{surfaces}. In some variants we assume that
the points lie on some fixed constant-degree algebraic variety, and in others
we consider arbitrary sets of points in 3-space.

The case of points and curves has been considered in several previous studies, starting with
Guth and Katz's work on points and lines~\cite{GK2}. Our results, which are based on
a recent work of Guth and Zahl~\cite{GZ} concerning surfaces that are doubly ruled by curves,
provide a grand generalization of all previous results. We reconstruct the bound
for points and lines, and improve, in certain signifcant ways, recent bounds involving
points and circles (in~\cite{SSZ}), and points and arbitrary constant-degree algebraic
curves (in~\cite{SSS}). While in these latter instances the bounds are not known (and
are strongly suspected not) to be tight, our bounds are, in a certain sense, the best
that can be obtained with this approach, given the current state of knowledge.

In the case of points and surfaces, the incidence graph between them
can contain large complete bipartite graphs, each involving points
on some curve and surfaces containing this curve (unlike earlier
studies, we do not rule out this possibility, which makes our
approach more general). Our bounds estimate the total size of the
\emph{vertex sets} in such a complete bipartite graph decomposition
of the incidence graph. In favorable cases, our bounds translate
into actual incidence bounds. Overall, here too our results provide
a ``grand generalization'' of most of the previous studies of
(special instances of) this problem.

As an application of our point-curve incidence bound, we consider
the problem of bounding the number of similar triangles spanned by a
set of $n$ points in $\reals^3$. We obtain the bound $O(n^{15/7})$,
which improves the bound of Agarwal et al.~\cite{AAPS}.

As applications of our point-surface incidence bounds, we consider the problems
of distinct and repeated distances determined by a set of $n$ points in $\reals^3$,
two of the most celebrated open problems in combinatorial geometry. We obtain new and
improved bounds for two special cases, one in which the points lie on some algebraic
variety of constant degree, and one involving incidences between pairs in $P_1\times P_2$,
where $P_1$ is contained in a variety and $P_2$ is arbitrary.
\end{abstract}


\section{Introduction} \label{sec:intro}

\subsection{The setups}

\paragraph{Incidences between points and curves in three dimensions.}
Let $P$ be a set of $m$ points and $\C$ a set of $n$ irreducible algebraic
curves of constant degree in $\reals^3$. We consider the problem of obtaining
sharp incidence bounds between the points of $P$ and the curves of $\C$.
This is a major topic in incidence geometry since the groundbreaking work of
Guth and Katz~\cite{GK2} on point-line incidences in $\reals^3$,
with many follow-up studies, some of which are reviewed below.
Building on the recent work of Guth and Zahl~\cite{GZ}, which bounds
the number of \emph{2-rich points} determined by a set of bounded-degree algebraic
curves in $\reals^3$ (i.e., points incident to at least two of the given curves),
we are able to generalize Guth and Katz's point-line incidence bound to a general
bound on the number of incidences between points and bounded-degree irreducible
algebraic curves that satisfy certain natural assumptions, discussed in detail below.

\paragraph{Incidences between points and surfaces in three dimensions.}
Let $P$ be a set of $m$ points, and $S$ a set of $n$ two-dimensonal algebraic
varieties of constant maximum degree in $\reals^3$. Here too we impose certain
natural assumptions on the surfaces in $S$, discussed in detail below.

Let $G(P,S) \subseteq P\times S$ denote the \emph{incidence graph} of $P$ and $S$;
its edges connect all pairs $(p,\sigma)\in P\times S$ such that $p$ is incident to $\sigma$.
In general, $I(P,S):= |G(P,S)|$ might be as large as the maximum possible value $mn$,
by placing all the points of $P$ on a suitable curve, and make all the surfaces of $S$
contain that curve,\footnote{%
  This situation can arise in many instances, for example in the case of planes
  (where many of them can intersect in a common line), or spheres (where many can
  intersect in a common circle), but there are also many cases where this is impossible.
  In this latter situation, which we do not yet know how to characterize in a simple and general manner, our analysis
  becomes sharper---see below.}
in which case $G(P,S) = P\times S$. The bound that we are going to
obtain will of course acknowledge this possibility, and will in fact
bypass it altogether. Concretely, rather than bounding $I(P,S)$, our
basic approach will represent $G(P,S)$ as a union of complete
bipartite subgraphs ${\displaystyle \bigcup_{\gamma\in\Gamma_0}
\big( P_\gamma\times S_\gamma\big)}$, and of a ``leftover'' subgraph
$G_0(P,S)$ (which, in certain cases, might be empty), and derive an
upper bound for the overall size of their \emph{vertex sets},
namely, a bound on
$$
J(P,S) := \sum_{\gamma\in\Gamma_0} \big( |P_\gamma| + |S_\gamma| \big) ,
$$
where the decomposition is over a set $\Gamma_0$ of constant-degree
algebraic curves $\gamma$ so that, for each $\gamma\in\Gamma_0$,
$P_\gamma = P\cap \gamma$ and $S_\gamma$ is the set of the surfaces
of $S$ that contain $\gamma$. (In some cases we will derive
different bounds on $\sum_\gamma |P_\gamma|$ and on $\sum_{\gamma}
|S_\gamma|$.) For the residual subgraph $G_0(P,S)$, we derive a
sharp bound on the actual number of incidences that it encodes
(namely, the number of its edges). This generalizes previous results
in which one had to require that $G(P,S)$ does not contain some
fixed-size complete bipartite graph, or (only for spheres or planes)
that the surfaces in $S$ be ``non-degenerate'' (\cite{ApS2,ET}; see
below).

\paragraph{Incidences between points on a variety and surfaces.}
An interesting special case is where $P$ is contained in some
two-dimensional algebraic variety (surface) $V$ of constant degree.
Besides being, as we believe, a problem of independent
interest, it arises as a key subproblem in our analysis of the
general case discussed above.

We assume that the surfaces of $S$ are taken from an
\emph{$s$-dimensional family of surfaces}, meaning that each of them
can be represented by a constant number of real parameters (e.g., by
the coefficients of the monomials of the polynomial whose zero set
is the surface), so that, in this parametric space, the points
representing the surfaces of $S$ lie in an $s$-dimensional algebraic
variety $\F$ of some constant degree (to which we refer as the
``complexity'' of $\F$). This assumption, which holds in practically
all applications, extends, in an obvious manner, to
lower-dimensional varieties (e.g., curves) and to higher dimensions;
see Sharir and Zahl~\cite{SZ} for a more thorough study of this
notion.

\subsection{Background}

\paragraph{Points and curves, the planar case.}
The case of incidences between points and curves has a rich history,
starting with the aforementioned case of points and lines in the
plane~\cite{CEGSW,Sze,ST}, where the worst-case tight bound on the
number of incidences is $\Theta(m^{2/3}n^{2/3}+m+n)$, where $m$ is
the number of points and $n$ is the number of lines. Still in the
plane, Pach and Sharir~\cite{PS} extended this bound to incidence
bounds between points and curves with $k$ \emph{degrees of freedom}.
These are curves with the property that, for each set of $k$ points,
there are only $\mu=O(1)$ curves that pass through all of them, and
each pair of curves intersect in at most $\mu$ points; $\mu$ is
called the \emph{multiplicity} (of the degrees of freedom).
\begin{theorem}[Pach and Sharir \cite{PS}] \label{th:PS}
Let $P$ be a set of $m$ points in $\reals^2$ and let $\C$ be a set
of $n$ bounded-degree algebraic curves in $\reals^2$ with $k$
degrees of freedom and with multiplicity $\mu$. Then
$$
I(P,\C) = O\left(m^{\frac{k}{2k-1}}n^{\frac{2k-2}{2k-1}}+m+n\right)
,
$$
where the constant of proportionality depends on $k$ and $\mu$.
\end{theorem}

\noindent{\bf Remark.} The result of Pach and Sharir holds for more
general families of curves, not necessarily algebraic, but, since
algebraicity will be assumed in higher dimensions, we assume it also
in the plane.

Except for the case $k=2$ (lines have two degrees of freedom), the
bound is not known, and strongly suspected not to be tight in the
worst case. Indeed, in a series of papers during the
2000's~\cite{ANPPSS, ArS, MT}, an improved bound has been obtained
for incidences with circles, parabolas, or other families of curves
with certain properties (see~\cite{ANPPSS} for the precise
formulation). Specifically, for a set $P$ of $m$ points and a set
$\C$ of $n$ circles, or parabolas, or similar curves~\cite{ANPPSS},
we have
\begin {equation}
\label {eq:impbo} I(P,\C) =
O(m^{2/3}n^{2/3}+m^{6/11}n^{9/11}\log^{2/11}(m^3/n)+m+n).
\end {equation}
Some further (slightly) improved bounds, over the bound in
Theorem~\ref{th:PS}, for more general families of curves in the
plane have been obtained by Chan~\cite{Chan, Chan2} and by
Bien~\cite{Bien}. They are, however, considerably weaker than the
bound in (\ref{eq:impbo}).

Recently, Sharir and Zahl~\cite{SZ} have considered general families
of constant-degree algebraic curves in the plane that belong to an
\emph{$s$-dimensional family of curves}. Similarly to the case of
surfaces, discussed above, this means that each curve in that family
can be represented by a constant number of real parameters, so that,
in this parametric space, the points representing the curves lie in
an $s$-dimensional algebraic variety $\F$ of some constant degree
(to which we refer, as above, as the ``complexity'' of $\F$).
See~\cite{SZ} for more details.
\begin{theorem}[Sharir and Zahl~\cite{SZ}]\label{incPtCu}
Let $\C$ be a set of $n$ algebraic plane curves that belong to
an $s$-dimensional family $\F$ of curves of maximum constant degree
$E$, no two of which share a common irreducible component, and let
$P$ be a set of $m$ points in the plane. Then, for any $\eps>0$,
the number $I(P,\C)$ of incidences between the points of
$P$ and the curves of $\C$ satisfies
\begin{equation*} \label{newIncBd}
 I(P,\C) = O\Big(m^{\frac{2s}{5s-4}} n^{\frac{5s-6}{5s-4}+\eps} + m^{2/3}n^{2/3} + m + n\Big) ,
\end{equation*}
where the constant of proportionality depends on $\eps$, $s$, $E$,
and the complexity of the family $\F$.
\end{theorem}
Except for the factor $O(n^\eps)$, this is a significant improvement
over the bound in Theorem~\ref{th:PS} (for $s\ge 3$), in cases where
the assumptions in Theorem~\ref{incPtCu} imply (as they often do)
that $\C$ has $k=s$ degrees of freedom. Concretely, when $k=s$, we
obtain an improvement, except for the factor $n^\eps$, for the entire
``meaningful'' range $n^{1/s} \le m \le n^2$, in which the bound is
superlinear. The factor $n^\eps$ makes the bound in \cite{SZ}
slightly weaker only when $m$ is close to the lower end $n^{1/s}$ of
that range. Note also that for circles (where $s=3$), the bound in
Theorem~\ref{incPtCu} nearly coincides with the slightly more
refined bound (\ref{eq:impbo}).

\paragraph{Incidences with curves in three dimensions.}
The seminal work of Guth and Katz~\cite{GK2} establishes the sharper
bound $O(m^{1/2}n^{3/4} + m^{2/3}n^{1/3}q^{1/3} + m + n)$ on the
number of incidences between $m$ points and $n$ lines in $\reals^3$,
provided that no plane contains more than $q$ of the given lines.
This has lead to many recent works on incidences between points and
lines or other curves in three and higher dimensions;
see~\cite{CPS,GZ,SSS,SSZ,SS4d,SS4dv} for a sample of these results.
Most relevant to our present study are the works of Sharir, Sheffer,
and Solomon~\cite{SSS} on incidences between points and curves in
any dimension, the work of Sharir, Sheffer, and Zahl~\cite{SSZ} on
incidences between points and circles in three dimensions, and the
work of Sharir and Solomon~\cite{SS4d} on incidences between points
and lines in four dimensions, as well as several other studies of
point-line incidences by the authors~\cite{SS3d, SS4dv}.

Of particular significance is the recent work of Guth and
Zahl~\cite{GZ} on the number of 2-rich points in a collection of
curves, namely, points incident to at least two of the given curves.
For the case of lines, Guth and Katz~\cite{GK2} have shown that the
number of such points is $O(n^{3/2})$, when no plane or regulus
contains more than $O(n^{1/2})$ lines. Guth and Zahl obtain the same
asymptotic bound for general algebraic curves, under analogous (but
stricter) restrictive assumptions.

The new bounds that we will derive require the extension to three dimensions of the notions
of having $k$ degrees of freedom and of being an $s$-dimensional family of curves.
The definitions of these concepts, as given above for the planar case, extend,
more or less verbatim, to three (or higher) dimensions, but, even in typical situations, these two concepts do not coincide anymore. For example, lines in three dimensions have two degrees of freedom, but they form a
$4$-dimensional family of curves (this is the number of parameters needed to
specify a line in $\reals^3$).

\paragraph{Points and surfaces.}
Many of the earlier works on point-surface incidences have only considered special
classes of surfaces, most notably planes and spheres (see below). The case of more
general surfaces has barely been considered, till the recent work of Zahl~\cite{Za},
who has studied the general case of incidences between $m$ points and $n$
bounded-degree algebraic surfaces in $\reals^3$ that have \emph{$k$ degrees of freedom}.
More precisely, in analogy with the case of curves, one needs to assume that for any
$k$ points there are at most $\mu=O(1)$ of the given surfaces that pass through all of them.
Zahl's bound is $O(m^{\frac{2k}{3k-1}}n^{\frac{3k-3}{3k-1}}+m+n)$, with the constant
of proportionality depending on $k$, $\mu$, and the maximum degree of the surfaces.

By B\'ezout's theorem, if we require every triple of the given surfaces to have finite
intersection, the number of intersection points would be at most $E^3$, where $E$ is
the degree of the surfaces. In particular, $E^3+1$ points would then have at most two
of the given surfaces passing through all of them. In many instances, though, the
actual number of degrees of freedom can be shown to be much smaller.

Zahl's bound was later generalized by Basu and Sombra~\cite{BS14} to incidences between points
and bounded-degree hypersurfaces in $\reals^4$ satisfying certain analogous conditions.

\paragraph{Points and planes.}
Although we will not specifically address this special case, we refer the reader
to the earlier works on this problem, going back to Edelsbrunner,
Guibas and Sharir~\cite{EGS}. More recently, Apfelbaum and
Sharir~\cite{ApS} (see also Brass and Knauer~\cite{BK} and Elekes
and T{\'{o}}th~\cite{ET}) have shown that if the incidence graph, for a set $P$ of $m$ points and a set $H$ of $n$ planes, does not contain a copy of $K_{r,s}$, for constant parameters $r$
and $s$, then $I(P,H)=O(m^{3/4}n^{3/4}+m+n)$. In more generality,
Apfelbaum and Sharir~\cite{ApS} have shown that if $I=I(P,H)$ is
significantly larger than this bound, then $G(P,H)$ must contain a
large complete bipartite subgraph $P'\times H'$, such that
$|P'|\cdot |H'| = \Omega(I^2/(mn)) - O(m+n)$. Moreover, as also
shown in~\cite{ApS} (slightly improving a similar result of Brass
and Knauer~\cite{BK}), $G(P,H)$ can be expressed as the union of
complete bipartite graphs $P_i\times H_i$ so that $\sum_i
(|P_i|+|H_i|) = O(m^{3/4}n^{3/4}+m+n)$. (This is a specialization to
the case $d=3$ of a similar result of \cite{ApS,BK} in any dimension
$d$, and it concurs with the approach followed in this paper for
more general scenarios.) Recently, Solomon and Sharir~\cite{SS16}
improved this bound substantially when all the points of $P$ lie on
a constant-degree variety $V$.

\paragraph{Points and spheres.}
Earlier works on the special case of point-sphere incidences have considered the general setup,
where the points of $P$ are arbitrarily placed in $\reals^3$.  Initial partial results
go back to Chung~\cite{Chung} and to Clarkson et al.~\cite{CEGSW}, and continue
with the work of Aronov et al.~\cite{APST}.  Later, Agarwal et al.~\cite{AAPS}
have bounded the number of \emph{non-degenerate} spheres with respect to a given
point set; their bound was subsequently improved by Apfelbaum and Sharir~\cite{ApS2}.\footnote{%
  Given a finite point set $P\subset\reals^3$ and a constant $0<\eta<1$,
  a sphere $\sigma \subset \reals^3$ is called $\eta$-degenerate (with
  respect to $P$), if there exists a circle $c \subset \sigma$ such that
  $|c\cap P| \ge \eta|\sigma \cap P|$.}

The aforementioned recent work of Zahl~\cite{Za} can be applied in
the case of spheres if one assumes that no three, or any larger but
constant number, of the spheres intersect in a common circle. In
this case the family has $k=3$ degrees of freedom --- any three
points determine a unique circle that passes through all of them,
and, by assumption, only $O(1)$ spheres contain that circle. Zahl's
bound then becomes $O(m^{3/4}n^{3/4}+m+n)$. In particular, this
bound holds for congruent (unit) spheres (where three such spheres
can never contain a common circle). The case of incidences with unit
spheres have also been studied in Kaplan et al.~\cite{KMSS}, with
the same upper bound; see also~\cite{SS16a}.

If many spheres of the family can intersect in a common circle, the
bound does no longer hold. The only earlier work that handled this
situation is by Apfelbaum and Sharir~\cite{ApS}, where it was
assumed that the given spheres are non-degenerate. In this case the
bound obtained in~\cite{ApS} is $O(m^{8/11}n^{9/11} + m + n)$.
Interestingly, this is also the bound that Zahl's result would have
yielded if the sphere had $k=4$ degrees of freedom, which however
they only ``almost have'': four generic points determine a unique
sphere that passes through all of them, but four co-circular points
determine an infinity of such spheres.

\paragraph{Distinct and repeated distances in three dimensions.}
The case of spheres is of particular interest, because it arises, in
a standard and natural manner, in the analysis of \emph{distinct and
repeated distances} determined by $n$ points in three dimensions
(see Section~\ref{sec:dd3}, where we use these well-known reductions
in our analysis). After Guth and Katz's almost complete solution of
the number of distinct distances in the plane~\cite{GK2}, the
three-dimensional case has moved to the research forefront. The
prevailing conjecture is that the lower bound is $\Omega(n^{2/3})$
(the best possible in the worst case), but the current record, due
to Solymosi and Vu~\cite{SoVu}, is still far smaller\footnote{This
follows by substituting the new lower bound $\Omega(n/\log n)$ of
Guth and Katz for distinct distances in the plane, in the recursive
analysis of~\cite{SoVu}.} (close to $\Omega(n^{3/5})$), and the
problem seems much harder than its two-dimensional counterpart.
Obtaining lower bounds for distinct distances using circles or
spheres has in general been suboptimal when compared with more
effective methods (such as in~\cite{GK2}), but here we use it effectively to obtain new
lower bounds (larger than $\Omega(n^{2/3})$) when the points lie on
a variety of fixed degree.

The status of the case of repeated distances is also far from being satisfactory.
The planar case is ``stuck'' with the upper bound $O(n^{4/3})$ of
Spencer et al.~\cite{SST} from the 1980's. This bound also holds for points on
the 2-sphere, and there it is tight in the worst case (when the repeated distance
is $1$, say, and the radius of the sphere is $1/\sqrt{2}$)~\cite{EHP}, but it is
strongly believed that in the plane the correct bound is close to linear. In three
dimensions, the aforementioned bound of~\cite{KMSS,Za} immediately implies the upper
bound $O(n^{3/2})$ on the number of repreated distances (a slight improvement over
the earlier bound of Clarkson et al.~\cite{CEGSW}), but the best known lower bound
is only $\Omega(n^{4/3}\log\log n)$~\cite{Erd}.

\subsection{Our results} \label{sec:res}

\paragraph{Incidences with curves.}
We first consider the problem of incidences between points and
algebraic curves. Before we state our results, we discuss three
notions that are used in these statements. These are the notions of
\emph{$k$ degrees of freedom} (already mentioned above), of
\emph{constructibility}, and of \emph{surfaces infinitely ruled by
curves.}

\paragraph{$k$ degrees of freedom.}
Let $\C_0$ be an infinite family of irreducible algebraic curves of constant
degree $E$ in $\RR^3$. Formally, in complete analogy with the planar case, we say
that $\C_0$ has $k$ \emph{degrees of freedom} with \emph{multiplicity} $\mu$,
where $k$ and $\mu$ are constants, if (i) for every tuple of $k$ points in $\RR^3$
there are at most $\mu$ curves of $\C_0$ that are incident to all $k$ points, and
(ii) every pair of curves of $\C_0$ intersect in at most $\mu$ points. As in~\cite{PS},
the bounds that we derive depend more significantly on $k$ than on $\mu$---see below.

We remark that the notion of $k$ degrees of freedom gets more
involved for surfaces, and raises several annoying technical issues.
For example, how many points does it take to define, say, a sphere
(up to a fixed multiplicity)? As already observed earlier, four
\emph{generic} points do the job (they define a unique sphere
passing through all four of them), but four co-circular points do
not.

While it seems possible to come up with some sort of working
definition, we bypass this issue in this paper, by defining this
notion, for a family $\F$ of surfaces, only with respect to a given
surface $V$, by saying that $\F$ has \emph{$k$ degrees of freedom
with respect to $V$} if the family of the irreducible components of
the curves $\{\sigma\cap V \mid \sigma\in\F\}$, counted
\emph{without} multiplicity, has $k$ degrees of freedom, in the
sense just defined. In the case of spheres, for example, this
definition gives, as is easily checked, four degrees of freedom when
$V$ is neither a plane nor a sphere, but only three when $V$ is a
plane or a sphere.

\paragraph{Constructibility.}
In the statements of the following theorems, we also assume that $\C_0$
is a \emph{constructible} family of curves. This notion generalizes
the notion of being algebraic, and is discussed in detail
in Guth and Zahl~\cite{GZ}. Informally, a set $Y\subset \cplx^d$ is constructible
if it is a Boolean combination of algebraic sets. The formal definition goes as
follows (see, e.g., Harris~\cite[Lecture 3]{Harris}).
For $z\in\cplx$, define $v(0)=0$ and $v(z)=1$ for $z\ne 0$. Then $Y\subseteq \cplx^d$,
for some fixed $d$, is a \emph{constructible set} if there exist a finite set of polynomials
$f_j : \cplx^d \to \cplx$, for $j = 1,\ldots,J_Y$, and a subset $B_Y \subset \{0,1\}^{J_Y}$,
so that $x \in Y$ if and only if $\left(v(f_1(x)),\ldots, v(f_{J_Y}(x))\right) \in B_Y$.

When we apply this definition to a set of curves, we think of them as points in some
parametric (complex) $d$-space, where $d$ is the number of parameters needed to specify
a curve. When $J_Y=1$ we get all the algebraic hypersurfaces (that admit the implied
$d$-dimensional representation) and their complements. An $s$-dimensional family of curves,
for $s<d$, is obtained by taking $J_Y=d-s$ and $B_Y=\{0\}^{J_Y}$. In doing so, the curves that we obtain are complete intersections. Following Guth and Zahl (see also a comment to that effect in the appendix), this involves no loss of generality, because every curve is contained in a curve that is a complete intersection. In what follows, when we
talk about constructible sets, we implicitly assume that the ambient dimension $d$ is constant.

The constructible sets form a Boolean algebra. This means that finite unions and
intersections of constructible sets are constructible, and the complement of a
constructible set is constructible.
Another fundamental property of constructible sets is that, over $\cplx$, the projection
of a constructible set is constructible; this is known as Chevalley's theorem
(see Harris~\cite[Theorem 3.16]{Harris} and Guth and Zahl~\cite[Theorem 2.3]{GZ}).
If $Y$ is a constructible set, we define the \emph{complexity} of $Y$ to be $\min(\deg f_1 + \cdots + \deg f_{J_Y})$, where the minimum
is taken over all representations of $Y$, as described above. As just observed,
constructibility of a family $\C_0$ of curves extends the notion of $\C_0$ being
$s$-dimensional. One of the main motivations for using the notion of constructible
sets (rather than just $s$-dimensionality) is the fact, established by
Guth and Zahl~\cite[Proposition 3.3]{GZ}, that the set $\C_{3,E}$ of \emph{irreducible}
curves of degree at most $E$ in complex 3-dimensional space (either affine or projective)
is a constructible set of constant complexity that depends only on $E$. Moreover,
Theorem~\ref{th:ruledgz}, one of the central technical tools that we use in our analysis
(see below for its statement and proof), holds for constructible families of curves.


\smallskip

\paragraph{The connection between degrees of freedom and constructibility/dimensionality.}
Loosely speaking, in the plane the number of degrees of freedom and the dimensionality of a family of curves tend to be equal. In three dimensions the situation is different. This is because the constraint that a curve $\gamma$ passes through a point $p$ imposes two equations on the parameters defining $\gamma$. We therefore expect the number of degrees of freedom to be half the dimensionality. A few instances where this is indeed the case are:
(i) Lines in three dimensions have two degrees of freedom, and they form a
$4$-dimensional family of curves (this is the number of parameters needed to
specify a line in $\reals^3$).
(ii) Circles in three dimensions have three degrees of
freedom, and they form a $6$-dimensional family of curves (e.g., one needs three
parameters to specify the plane containing the circle, two additional parameters to
specify its center, and a sixth parameter for its radius). (iii) Ellipses have five
degrees of freedom, but they form an $8$-dimensional family of curves, as is easily checked.
(This discrepancy (for ellipses) is explained by noting that four points are not sufficient to define the ellipse because the first three determine the plane containing it, so the fourth point, if at all coplanar with the first three, only imposes one constraint on the parameters of the ellipse.)

\noindent{\bf Remark.} The definition of constructibility is given
over the complex field $\cplx$. This is in accordance with most of
the basic algebraic geometry tools, which have been developed over
the complex field. Some care has to be exercised when applying them
over the reals. For example, Theorem~\ref{th:ruledgz}, one of the
central technical tools that we use in our analysis, as well as the
results of Guth and Zahl~\cite{GZ}, apply over the complex field,
but not over the reals. On the other hand, when we apply the
partitioning method of~\cite{GK2} (as in the proofs of
Theorems~\ref{th:mainInc3}) and~\ref{th:incgen} or when we use
Theorem~\ref{incPtCu}, we (have to) work over the reals.

It is a fairly standard practice in algebraic geometry that handles
a real algebraic variety $V$, defined by real polynomials, by
considering its complex counterpart $V_\cplx$, namely the set of
complex points at which the polynomials defining $V$ vanish. The
rich toolbox that complex algebraic geometry has developed allows
one to derive various properties of $V_\cplx$, which, with some care, can
usually be transported back to the real variety $V$.

This issue arises time and again in this paper. Roughly speaking, we approach
it as follows. We apply the polynomial partitioning technique to the given
sets of points and of curves or surfaces, in the original real (affine) space,
as we should. Within the cells of the partitioning we then apply some
field-independent argument, based either on induction or on some ad-hoc
combinatorial argument. Then we need to treat points that lie on the zero set
of the partitioning polynomial. We can then switch to the complex field, when
it suits our purpose, noting that this step preserves all the real incidences;
at worst, it might add additional incidences involving the non-real portions
of the variety and of the curves or surfaces. Hence, the bounds that we obtain
for this case transport, more or less verbatim, to the real case too.

\paragraph{Surfaces infinitely ruled by curves.}
Back in three dimensions, a surface $V$ is (singly, doubly, or
infinitely) ruled by some family $\Gamma$ of curves of degree at
most $E$, if each point $p\in V$ is incident to (at least one, at
least two, or infinitely many) curves of $\Gamma$ that are fully
contained in $V$. The connection between ruled surface theory and
incidence geometry goes back to the pioneering work of Guth and
Katz~\cite{GK2} and shows up in many subsequent works. See Guth's
recent survey~\cite{Gut:surv} and recent book~\cite{Gut:book}, and
Koll\'ar~\cite{Kollar} for details.

In most of the previous works, only singly-ruled and doubly-ruled
surfaces have been considered. Looking at infinitely-ruled surfaces
adds a powerful ingredient to the toolbox, as will be demonstrated
in this paper.

We recall that the only surfaces that are infinitely ruled by lines
are planes (see, e.g., Fuchs and Tabachnikov~\cite[Corollary
16.2]{FT}), and that the only surfaces that are infinitely ruled
by circles are spheres and planes (see, e.g., Lubbes~\cite[Theorem
3]{Lub} and Schicho~\cite{Sch}; see also Skopenkov and
Krasauskas~\cite{SK} for recent work on \emph{celestials}, namely
surfaces \emph{doubly} ruled by circles, and Nilov and
Skopenkov~\cite{NS13}, proving that a surface that is ruled by a
line and a circle through each point is a quadric). It should be noted that,
in general, for this definition to make sense, it is important to
require that the degree $E$ of the ruling curves be much smaller
than $\deg(V)$. Otherwise, every variety $V$ is infinitely ruled by,
say, the curves $V\cap h$, for hyperplanes $h$, having the same
degree as $V$. A challenging open problem is to characterize all the
surfaces that are infinitely ruled by algebraic curves of degree at
most $E$ (or by certain special classes thereof). However, the
following result of Guth and Zahl provides a useful sufficient
condition for this property to hold.
\begin{theorem} [Guth and Zahl~\cite{GZ}] \label{th:gz}
Let $V$ be an irreducible surface, and suppose that it is
\emph{doubly ruled} by curves of degree at most $E$. Then $\deg(V) \le 100E^2$.
\end{theorem}
In particular, an irreducible surface that is infinitely ruled by
curves of degree at most $E$ is doubly ruled by these curves, so its
degree is at most $100E^2$. Therefore, if $V$ is irreducible of
degree $D$ larger than this bound, $V$ cannot be infinitely ruled by
curves of degree at most $E$. This leaves a gray zone, in which the
degree of $V$ is between $E$ and $100E^2$. We would like to
conjecture that in fact no irreducible variety with degree in this
range is infinitely ruled by degree-$E$ curves. Being unable to establish this conjecture, we leave it as a challenging open problem for further research.

Finally, we remark that the notion of surfaces infinitely ruled by curves also plays a crucial
role in one of our results on point-surface incidences (see Theorem~\ref{incth}).

\paragraph{Our results: points and curves.}

We can now state our main results on point-curve incidences.
\begin{theorem}[Curves in $\RR^3$] \label{th:mainInc3}
Let $P$ be a set of $m$ points and $\C$ a set of $n$ irreducible
algebraic curves of constant degree $E$, taken from a constructible family
$\C_0$, of constant complexity, with $k$ degrees of freedom (and some multiplicity
$\mu$) in $\RR^3$, such that no surface that is infinitely ruled by curves of $\C_0$
contains more than $q$ curves of $\C$, for a parameter $q<n$. Then
\begin{equation} \label{eq:mainInc3}
I(P,\C) = O\left(m^{\frac{k}{3k-2}}n^{\frac{3k-3}{3k-2}} +
m^{\frac{k}{2k-1}}n^{\frac{k-1}{2k-1}}q^{\frac{k-1}{2k-1}}+m+n\right),
\end{equation}
where the constant of proportionality depends on $k$, $\mu$, $E$, and
the complexity of the family $\C_0$.
\end{theorem}

\noindent{\bf Remarks. (1)} In certain favorable situations, such as
in the cases of lines or circles, discussed above, the surfaces that
are infinitely ruled by curves of $\C_0$ have a simple
characterization. In such cases the theorem has a stronger flavor,
as its assumption on the maximum number of curves on a surface has
to be made only for this concrete kind of surfaces. For example, as
already noted, for lines (resp., circles) we only need to require
that no plane (resp., no plane or sphere) contains more than $q$ of
the curves. In general, as mentioned, characterizing
infinitely-ruled surfaces by a specific family of curves is a
difficult task. Nevertheless, we can overcome this issue by
replacing the assumption in the theorem by a more restrictive one,
requiring that no surface that is infinitely ruled by curves of
degree at most $E$ contain more than $q$ curves of $\C$. By
Theorem~\ref{th:gz}, any infinitely ruled surface of this kind must
be of degree at most $100E^2$. Hence, an even simpler (albeit
weaker) formulation of the theorem is to require that no surface of
degree at most $100E^2$ contains more than $q$ curves of $\C$. This
can indeed be much weaker: In the case of circles, say, instead of
making this requirement only for planes and spheres, we now have to
make it for every surface of degree at most $400$.

\smallskip

\noindent{\bf (2)} In several recent works (see~\cite{Gut,SSS,SSZ}),
the assumption in the theorem is replaced by a much more restrictive
assumption, that no surface of degree at most $c_\eps$ contains more
than $q$ given curves, where $c_\eps$ is a constant that depends on
another prespecified parameter $\eps>0$ (where $\eps$ appears in the
exponents in the resulting incidence bound), and is typically very
large (and increases as $\eps$ becomes smaller). Getting rid of such
an $\eps$-dependent constant (and of the $\eps$ in the exponent) is a significant feature of
Theorem~\ref{th:mainInc3}.

\smallskip

\noindent{\bf (3)}
Theorem~\ref{th:mainInc3} generalizes the incidence bound of Guth and Katz~\cite{GK2},
obtained for the case of lines. In this case, lines have $k=2$ degrees of freedom,
they certainly form a constructible (in fact, a $4$-dimensional) family of curves,
and, as just noted, planes are the only surfaces in $\RR^3$ that are infinitely
ruled by lines. Thus, in this special case, both the assumptions and the bound in
Theorem~\ref{th:mainInc3} are identical to those in Guth and Katz~\cite{GK2}.
That is, if no plane contains more than $q$ input lines, the number of incidences is
$O(m^{1/2}n^{3/4}+m^{2/3}n^{1/3}q^{1/3}+m+n)$.

\paragraph{Improving the bound.}
The bound in Theorem~\ref{th:mainInc3} can be further improved, if
we also throw into the analysis the dimensionality $s$ of the family
$\C_0$. Actually, as will follow from the proof, the dimensionality
that will be used is only that of any subset of $\C_0$ whose members
are fully contained in some variety that is infinitely ruled by
curves of $\C_0$. As just noted, such a variety must be of constant
degree (at most $100E^2$, or smaller as in the cases of lines and
circles), and the additional constraint that the curves be contained
in the variety can typically be expected to reduce the
dimensionality of the family.

For example, if $\C_0$ is the collection of all circles in $\reals^3$, then, since
the only surfaces that are infinitely ruled by circles are spheres and planes, the
subfamily of all circles that are contained in some sphere or plane is only
$3$-dimensional (as opposed to the entire $\C_0$, which is $6$-dimensional).

We capture this setup by saying that $\C_0$ is a family of \emph{reduced dimension}
$s$ if, for each surface $V$ that is infinitely ruled by curves of $\C_0$, the subfamily
of the curves of $\C_0$ that are fully contained in $V$ is $s$-dimensional. In this
case we obtain the following variant of Theorem~\ref{th:mainInc3}.
\begin{theorem}[Curves in $\RR^3$] \label{th:imprInc3}
Let $P$ be a set of $m$ points and $\C$ a set of $n$ irreducible
algebraic curves of constant degree $E$, taken from a constructible family
$\C_0$ with $k$ degrees of freedom (and some multiplicity $\mu$) in
$\RR^3$, such that no surface that is infinitely ruled by curves of
$\C_0$ contains more than $q$ of the curves of $\C$, and assume further that
$\C_0$ is of reduced dimension $s$. Then
\begin{equation} \label{eq:imprInc3}
I(P,\C) = O\left(m^{\frac{k}{3k-2}}n^{\frac{3k-3}{3k-2}}
\right) + O_\eps\left(
m^{2/3}n^{1/3}q^{1/3} +
m^{\frac{2s}{5s-4}}n^{\frac{3s-4}{5s-4}}q^{\frac{2s-2}{5s-4}+\eps}+m+n\right),
\end{equation}
for any $\eps>0$, where the first constant of proportionality depends on
$k$, $\mu$, $s$, $E$, and the maximum complexity of any subfamily
of $\C_0$ consisting of curves that are fully contained in some surface
that is infinitely ruled by curves of $\C_0$, and the second constant also depends on $\eps$.
\end {theorem}

\noindent{\bf Remarks. (1)}
Theorem~\ref{th:imprInc3} is an improvement of Theorem~\ref{th:mainInc3} when
$s\le k$ and $m>n^{1/k}$, in cases where $q$ is sufficiently large so as to make
the second term in (\ref{eq:mainInc3}) dominate the first term; for smaller values
of $m$ the bound is always linear. This is true except for the term $q^\eps$, which
affects the bound only when $m$ is very close to $n^{1/k}$ (when $s=k$). When $s>k$
we get a threshold exponent $\beta = \frac{5s-4k-2}{ks-4k+2s}$ (which becomes $1/k$
when $s=k$), so that the bound in Theorem~\ref{th:imprInc3} is stronger (resp., weaker)
than the bound in Theorem~\ref{th:mainInc3} when $m > n^\beta$ (resp., $m < n^\beta$),
again, up to the extra factor $q^\eps$.

\smallskip

\noindent{\bf (2)}
The bounds in Theorems~\ref{th:mainInc3} and~\ref{th:imprInc3} improve, in three
dimensions, the recent result of Sharir, Sheffer, and Solomon~\cite{SSS}, in three
significant ways:
{\bf(i)} The leading terms in both bounds are essentially the same, but
our bound is sharper, in that it does not include the factor $O(n^\eps)$ appearing
in \cite{SSS}. {\bf (ii)} The assumption here, concerning the number of curves on a
low-degree surface, is much weaker than the one made in~\cite{SSS}, where
it was required that no surface of some (constant but potentially very large)
degree $c_\eps$, that depends on $\eps$, contains more than $q$ curves of $\C$.
(See also Remark (2) following Theorem~\ref{th:mainInc3}.)
{\bf (iii)} The two variants of the non-leading terms here are significantly smaller
than those in~\cite{SSS}, and, in a certain sense (that will be elaborated
following the proof of Theorem~\ref{th:imprInc3}) are best possible.

\smallskip
\paragraph{Point-circle incidences in $\reals^3$.}

Theorem~\ref{th:imprInc3} yields a new bound for the case of incidences between points and circles
in $\reals^3$, which improves over the previous bound of Sharir,
Sheffer, and Zahl~\cite{SSZ}. Specifically, we have:
\begin{theorem} \label{ptcirc}
Let $P$ be a set of $m$ points and $\C$ a set of $n$ circles in $\reals^3$, so that
no plane or sphere contains more than $q$ circles of $\C$. Then
$$
I(P,\C) = O\left( m^{3/7}n^{6/7} + m^{2/3}n^{1/3}q^{1/3} +
m^{6/11}n^{5/11}q^{4/11}\log^{2/11}(m^3/q) + m + n \right) .
$$
\end{theorem}
Here too we have the three improvements noted in Remark (2) above. In particular, in the sense of part (iii)
of that remark, the new bound is ``best possible'' with respect to the best known bound (\ref{eq:impbo}) for
the planar or spherical cases. See Section~\ref{se:sim} for details.
Theorem~\ref{ptcirc} has an interesting application to the problem of bounding the number of similar triangles spanned by a set of $n$ points in $\reals^3$. It yields the bound $O(n^{15/7})$, which improves the bound of Agarwal et al.~\cite{AAPS}. See Section~\ref{se:sim} for details.
\paragraph{Incidence graph decomposition, for points on a variety and surfaces.}
Our first main result on point-surface incidences deals with the
special case where the points of $P$ lie on some algebraic variety
$V$ of constant degree. Besides being of independent interest, this
is a major ingredient of the analysis for the general case of an
arbitrary set of points in $\reals^3$ and surfaces.

In the statements of the following theorems we assume that the set
$S$ of the given surfaces is taken from some infinite family $\F$
that either has $k$ degrees of freedom with respect to $V$ (with
some multiplicity $\mu$), as defined earlier, for suitable constant parameters
$k$ (and $\mu$), or is of \emph{reduced
dimension} $s$ with respect to $V$, for some constant parameter $s$, meaning that the
family $\Gamma:= \{\sigma\cap V \mid \sigma\in F\}$ is an
$s$-dimensional family of curves (this is reminiscent of the notion
of reduced dimension defined above for curves).
\begin{theorem} \label{main}
Let $P$ be a set of $m$ points on some algebraic surface $V$ of
constant degree $D$ in $\reals^3$, and let $S$ be a set of $n$
algebraic surfaces in $\reals^3$ of maximum constant degree $E$,
taken from some family $\F$ of surfaces, which either has $k$
degrees of freedom with respect to $V$ (with some multiplicity
$\mu$), or is of reduced dimension $s$ with respect to $V$, for some
constant parameters $k$ (and $\mu$) or $s$. We also assume that the
surfaces in $S$ do not share any common irreducible component (which
certainly holds when they are irreducible). Then the incidence graph
$G(P,S)$ can be decomposed as
\begin{equation} \label{gps}
G(P,S) = \bigcup_\gamma (P_\gamma\times S_\gamma) ,
\end{equation}
where the union is over all irreducible components of curves
$\gamma$ of the form $\sigma\cap V$, for $\sigma\in S$, and, for
each such $\gamma$, $P_\gamma = P\cap \gamma$ and $S_\gamma$ is the
set of surfaces in $S$ that contain $\gamma$.

If $\F$ has $k$ degrees of freedom then
\begin{equation} \label{eq:maink}
\sum_\gamma |P_\gamma| = O\Big(m^{\frac{k}{2k-1}}
n^{\frac{2k-2}{2k-1}} + m + n\Big) ,
\end{equation}
and if $\F$ is $s$-dimensional then we have, for any $\eps>0$,
\begin{equation} \label{eq:mains}
\sum_\gamma |P_\gamma| = O\Big(m^{\frac{2s}{5s-4}}
n^{\frac{5s-6}{5s-4}+\eps} + m^{2/3}n^{2/3} + m + n\Big) ,
\end{equation}
where the constants of proportionality depends on $D$, $E$, and the
complexity of the family $\F$, and either on $k$ and $\mu$ in the
former case, or on $\eps$ and $s$ in the latter case.

Moreover, in both cases we have $\sum_\gamma |S_\gamma| = O(n)$, where the
constant of proportionality depends on $D$ and $E$.

\end{theorem}
\smallskip

\noindent {\bf Remark.} A major feature of this result is that it
does not impose any restrictions on the incidence graph, such as
requiring it not to contain some fixed complete bipartite graph
$K_{r,r}$, for $r$ a constant, as is done in the preceding
studies~\cite{BS14,KMSS,Za}. We re-iterate that, to allow for the
existence of large complete bipartite graphs, the bounds in
(\ref{eq:maink}) and (\ref{eq:mains}), as well as the bound
$\sum_\gamma |S_\gamma|=O(n)$, are not on the number of incidences
(that is, on the number of edges in $G(P,S)$, which could be as high
as $mn$) but on the overall size of the vertex sets of the subgraphs
in the complete bipartite graph decomposition of $G(P,S)$. This
would lead to the same asymptotic bound on $|G(P,S)|$ itself, if one
assumes that this graph does not contain $K_{r,r}$ as a subgraph,
for a constant $r$.

This kind of compact representation of incidences has already been used in the previous
studies of Brass and Knauer~\cite{BK}, Apfelbaum and Sharir~\cite{ApS}, and our recent
works~\cite{SS16, SS16a}, albeit only for the special cases of planes or spheres.

\noindent {\bf Remark.} Another way of bypassing the possible
presence of large complete bipartite graphs in $G(P,S)$, used in
several earlier works~\cite{AAPS,ApS,ET}, is to assume that the
surfaces in $S$ are \emph{non-degenerate}. These studies, already
mentioned earlier, only considered the cases of planes and spheres
(or of hyperplanes and spheres in higher dimensions)~\cite{AAPS,
ET}. For spheres, for example, this means that no more than some
fixed fraction of the points of $P$ on any given sphere can be
cocircular. Although large complete bipartite graphs can exist in
$G(P,S)$ in this case, the non-degeneracy assumption allows us to
control, in a sharp form, the number of incidences (and shows that
the resulting complete bipartite graphs are not so large after all).
It would be interesting (and, as we believe, doable) to extend our
analysis to the case of (suitably defined) more general
non-degenerate surfaces. These remarks also apply to the general
case (involving points anywhere in $\reals^3$), given in
Theorem~\ref{th:incgen} below.

\paragraph{A mixed incidence bound (for points on most varieties and general surfaces).}
Our second result is an improvement of Theorem~\ref{main}, still for
the case where the points of $P$ lie on some algebraic variety $V$
of constant degree, where we now also assume that $V$ is not
\emph{infinitely ruled} by the (irreducible components of the)
intersection curves of pairs of members of the given family $\F$ of surfaces.
In this case we obtain an improved, ``mixed'' bound, in which
$G(P,S)$ can be split into two subgraphs, $G_0(P,S)$ and $G_1(P,S)$,
where the bound in~(\ref{eq:maink}) or in~(\ref{eq:mains}) now holds
for $|G_0(P,S)|$, i.e., for the actual number of incidences that it
represents, and where $G_1(P,S)$ admits a complete bipartite graph
decomposition, as
above, for which the sum of the vertex sets is only\footnote{%
  In fact, many ``bad'' things must happen for $G_1(P,S)$ to be nontrivial,
  and in many situations one would expect $G_1(P,S)$ to be empty; see below.}
$O(m+n)$. The actual bound is slightly sharper---see below.

Specializing the theorem to the case of spheres, as is done later on
(in Section~\ref{sec:dd3}), leads to interesting implications to
distinct and repeated distances in three dimensions.

\begin{theorem} \label{incth}
Let $P$ be a set of $m$ points on some irreducible algebraic surface
$V$ of constant degree $D$ in $\reals^3$, and let $S$ be a set of
$n$ algebraic surfaces in $\reals^3$ of constant degree $E$, which
do not share any common irreducible component, taken from some
infinite constructible family $\F$ of surfaces that either has $k$
degrees of freedom with respect to $V$ (with some multiplicity
$\mu$) or is $s$-dimensional with respect to $V$, for some constant
parameters $k$ (and $\mu$) or $s$. Assume further that $V$ is not
infinitely ruled by the family $\C_0$ of the irreducible components
of the intersection curves of pairs of surfaces\footnote{A stricter
assumption is that $V$ is not infinitely ruled by algebraic curves
of degree at most $E^2$, which will hold if we assume that each
irreducible component of $V$ has degree larger than $100E^4$.} in
$\F$. Then the incidence graph $G(P,S)$ can be decomposed as
\begin{equation} \label{gps1}
G(P,S) = G_0(P,S) \cup \bigcup_\gamma (P_\gamma\times S_\gamma) ,
\end{equation}
where the union is over all irreducible curves $\gamma$ contained in
(one-dimensional) intersections of the form $\sigma\cap \sigma'\cap V$,
for $\sigma\ne\sigma'\in S$, and, for each such $\gamma$,
$P_\gamma \subseteq P\cap \gamma$ (for some points on some curves, their incident pairs are moved to, and counted in $G_0(P,S)$), and $S_\gamma$ is the set (of size at least two)
of surfaces in $S$ that contain $\gamma$.

Moreover, if $\F$ has $k$ degrees of freedom with respect to $V$
(with some multiplicity $\mu$) then
\begin{equation} \label{eq:main2k}
|G_0(P,S)| = O\Big(m^{\frac{k}{2k-1}} n^{\frac{2k-2}{2k-1}} + m + n\Big) ,
\end{equation}
and if $\F$ is $s$-dimensional with respect to $V$ then, for any
$\eps>0$,
\begin{equation} \label{eq:main2s}
|G_0(P,S)| = O\Big(m^{\frac{2s}{5s-4}} n^{\frac{5s-6}{5s-4}+\eps} +
m^{2/3}n^{2/3} + m + n\Big) ,
\end{equation}
where the constants of proportionality depends on $D$, $E$, and the complexity of the family $\F$,
and either on $k$ and $\mu$ in the former case, or on $\eps$ and $s$ in the latter case.
In either case we also have
$$
\sum_\gamma |P_\gamma|=O(m), \quad\text{and}\quad \sum_\gamma
|S_\gamma|=O(n) \nonumber ,
$$
where the constants of proportionality depend on $D$, $E$, and
the complexity of the family $\F$, and either on $k$ (and $\mu$) in
the former case, or on $\eps$ and $s$ in the latter case.
\end{theorem}


\smallskip

\noindent{\bf Remarks. (1)} As already alluded to, we note that,
typically, one would expect the complete bipartite decomposition
part of (\ref{gps1}) to be empty or trivial. To really be
significant, (a) many surfaces of $S$ would have to intersect in a
common curve, and, in cases where the multiplicity of these curves
is not that large, (b) many curves of this kind would have to be
fully contained in $V$. Thus, in many cases, in which (a) and (b) do
not hold, the bounds in (\ref{eq:main2k}) or in (\ref{eq:main2s}) in
Theorem~\ref{incth} are for the overall number of incidences. Note
also that both Theorem~\ref{main} and Theorem~\ref{incth} yield a
decomposition of (the whole or a portion of) $G(P,S)$ into complete
bipartite subgraphs. The major difference is that the bound
$\sum_\gamma |P_\gamma|$ on the overall $P$-vertex sets size of
these graphs is (relatively) large in Theorem~\ref{main}, but it is
only \emph{linear} in $m$ and $n$ (if at all nonzero) in
Theorem~\ref{incth}. (The bound on $\sum_{\gamma} |S_\gamma|$
remains $O(n)$ in both cases.)

\smallskip

\noindent{\bf (2)} We note that if $V$ is infinitely ruled by our
curves the results break down. For a simple example, take $m$ points
and $N$ lines in the plane which form $\Theta(m^{2/3}N^{2/3})$
incidences between them. Now pick any surface $V$ in $\reals^3$, say
the paraboloid $z=x^2+y^2$ for specificity, and lift up each of the
$N$ lines to a vertical parabola on $V$. Clearly, $V$ is infinitely
ruled by such parabolas, and we get a system of $m$ points and $n$
parabolas with $\Theta(m^{2/3}N^{2/3})$ incidences between them. It
is also easy to turn this construction into a point-surface
incidence structure, in which $\sum_\gamma |P_\gamma|$ is equal to
this bound, which is larger than the lower bound $O(m+N)$ asserted
in the theorem. The line $y=ax+b$ in the plane is lifted to the
parabola $\gamma_{a,b} =\{(x,y,z)\in \reals^3: y=ax+b, z=x^2+y^2\}$
contained in the paraboloid $V$. Define a family $S$ of quadratic
surfaces parameterized by $a,b,c_0,c_1,c_2 \in \reals$ by
$S_{a,b,c_0,c_1,c_2} :=\{(x,y,z)\in \reals^3 \mid (z-x^2-y^2) +
(y-ax-b)(c_0+c_1x+c_2y) = 0\}$. For any $c_0,c_1,c_2 \in \reals$,
the quadric $S_{a,b,c_0,c_1,c_2}$ contains the parabola
$\gamma_{a,b}$, i.e., many surfaces in $S$ intersect in a common
parabola.

\paragraph{Incidences between points on a variety and spheres.}
A particular case of interest is when $S$ is a set of spheres. The
intersection curves of spheres are circles, and, as already noted,
the only surfaces that are infinitely ruled by circles are spheres
and planes. Hence, to apply Theorem~\ref{incth}, we need to assume
that the constant-degree surface $V$ that contains the points of $P$
has no planar or spherical components, thereby ensuring that $V$ is
not infinitely ruled by circles. Clearly, as already noted, spheres
in $\reals^3$ have four degrees of freedom with respect to any
constant-degree variety with no planar or spherical components, and
they form a four-dimensional family of surfaces, with respect to any
such variety (and also in general). We can therefore apply
Theorem~\ref{incth}, with $s=4$, and conclude:
\begin{theorem} \label{main:sph}
Let $P$ be a set of $m$ points on some algebraic surface $V$ of
constant degree $D$ in $\reals^3$, which has no linear or spherical
components, and let $S$ be a set of $n$ spheres, of arbitrary radii,
in $\reals^3$. The incidence graph $G(P,S)$ can be decomposed as
\begin{equation} \label{gps-sph}
G(P,S) = G_0(P,S) \cup \bigcup_{\gamma\in\Gamma} (P_\gamma\times S_\gamma) ,
\end{equation}
where $\Gamma$ is the set of circles that are contained in $V$ and in at least two spheres of $S$,
and such that, for each $\gamma\in\Gamma$, $P_\gamma = P\cap\gamma$ and $S_\gamma$ is the
set of all spheres in $S$ that contain $\gamma$. We have
\begin{align} \label{eq:mainsph}
|G_0(P,S)| & = O\left( m^{1/2}n^{7/8+\eps} + m^{2/3}n^{2/3} + m + n \right) , \\
\sum_\gamma |P_\gamma| & = O(m) , \quad\text{and}\quad \sum_\gamma |S_\gamma| = O(n) ,
\nonumber ,
\end{align}
for any $\eps>0$, where the constant of proportionality depends on $D$ and $\eps$.
\end{theorem}

\noindent{\bf Remark.} Since $V$ does not contain a planar or
spherical component, the number of circles in $\Gamma$ is $O(D^2)$,
as follows by Guth and Zahl~\cite{GZ}. That is, the union in
(\ref{eq:mainsph}) is only over a constant number of circles. On the
other hand, there might also be incidence edges contained in
complete bipartite graphs corresponding to circles that are not
contained in $V$, whose number might be quite large. These
incidences are recorded in $G_0(P,S)$ and their number is bounded in
(\ref{eq:mainsph}).

Zahl's assumption that $G(P,S)$ does not contain $K_{r,3}$, for
some (arbitrary) constant $r$ (that is, by assuming that every
triple of spheres intersect in at most $r$ points of $P$), leads to the bound $I(P,S) = O(m^{3/4}n^{3/4}+m+n)$; our bound is better for
$m>n^{1/2}$ (ignoring the $n^{\eps}$ factor in our bound). Except
for this rather restrictive assumption, Zahl's result is more
general, as it does not require the points to lie on a
constant-degree variety.

We also note that if we assume that $G(P,S)$ does not contain any
$K_{r,r}$, for $r>3$ a constant, the bound in the second part of (\ref{eq:mainsph}) becomes a bound on the number of incidences, so, under this somewhat weaker assumption (than that of Zahl), we
improve Zahl's bound for points on a variety and for $m>n^{1/2}$.

The bound in (\ref{eq:mainsph}) further improves when either (i) the centers
of the spheres of $S$ lie on $V$ (or on some other constant-degree variety),
or (ii) the spheres of $S$ have the same radius. In both cases, $S$ is only
three-dimensional, so the bound improves to
\begin{equation} \label{eq:mainsph2}
|G_0(P,S)| = O\left( m^{6/11}n^{9/11+\eps} + m^{2/3}n^{2/3} + m + n \right) ,
\end{equation}
for any $\eps>0$. When both conditions hold---the spheres are
congruent and their centers lie on $V$---$S$ is only two-dimensional
with respect to $V$, and the bound improves still further to
\begin{equation*}
|G_0(P,S)| = O\left( m^{2/3}n^{2/3 + \eps} + m + n \right) .
\end{equation*}
Using a slightly refined machinery, developed in a companion
paper~\cite{SS16a}, the latter bound can be actually improved further to
\begin{equation} \label{eq:mainsph3}
O(m^{2/3}n^{2/3}+m+n).
\end{equation}

\paragraph{\bf Applications of Theorem~\ref{main:sph} and (\ref{eq:mainsph2}), (\ref{eq:mainsph3}):}

\smallskip

\noindent{\bf Distinct distances.}
As already mentioned, and as will be detailed in the proofs of the following results,
the new bounds on point-sphere incidences have immediate applications to the study of
distinct and repeated distances determined by a set of $n$ points in $\reals^3$, when
the points (or a subset thereof---see below) lie on some fixed-degree algebraic variety. Specifically, for distinct
distances, we have the following results.
\begin{theorem} \label{dd3}
(a) Let $P$ be a set of $n$ points on an algebraic surface $V$ of
constant degree $D$ in $\reals^3$, with no linear or spherical
components. Then the number of distinct distances determined by $P$
is $\Omega(n^{7/9-\eps})$, for any $\eps>0$, where the constant of proportionality depends on $D$ and $\eps$.

\smallskip

\noindent (b) Let $P_1$ be a set of $m$ points on a surface $V$ as
in (a), and let $P_2$ be a set of $n$ arbitrary points in
$\reals^3$. Then the number of distinct distances determined by
pairs of points in $P_1\times P_2$ is
$$
\Omega \left( \min\left\{ m^{4/7-\eps}n^{1/7-\eps},\; m^{1/2}n^{1/2},\; m \right\} \right) ,
$$
for any $\eps>0$, where the constant of proportionality depends on $D$ and $\eps$.
\end{theorem}

\smallskip

\noindent{\bf Remark.} In a recent work~\cite{SS16a}, we have
obtained slightly improved bounds, without the $\eps$ in the
exponents, using a more refined space decomposition technique, which
can be applied for arrangements of spheres.

While we believe that the bounds in the theorem are not tight, we
note that the bounds in both (a) and (b) (with, say, $m=n$) are
significantly larger than the conjectured best-possible lower bound
$\Omega(n^{2/3})$ for arbitrary point sets in $\reals^3$.

\smallskip

\noindent{\bf Repeated distances.} As another application, we bound
the number of unit (or repeated) distances involving points on a
surface $V$, as above.

\begin{theorem} \label{und}
(a) Let $P$ be a set of $n$ points on some algebraic surface $V$ of
constant degree $D$ in $\reals^3$, which does not contain any planar
or spherical components.
Then $P$ determines $O(n^{4/3})$ unit distances, where the constant of proportionality depends on $D$.

\smallskip

\noindent (b) Let $P_1$ be a set of $m$ points on a surface $V$ as
in (a), and let $P_2$ be a set of $n$ arbitrary points in
$\reals^3$. Then the number of unit distances determined by pairs of
points in $P_1\times P_2$ is
$$
O\left( m^{6/11}n^{9/11+\eps} + m^{2/3}n^{2/3} + m + n \right) ,
$$
for any $\eps>0$, where the constant of proportionality depends on $D$ and $\eps$.
\end{theorem}
In part (a) we extend, to the case of general constant-degree algebraic surfaces, the known
bound $O(n^{4/3})$, which is worst-case tight when $V$ is a sphere~\cite{EHP}. Part (b) gives
(say, for the case $m=n$) an intermediate bound between $O(n^{4/3})$ and the best known upper
bound $O(n^{3/2})$ for a arbitrary set of points in $\reals^3$~\cite{KMSS,Za}.

Another thing to notice is that, for distinct distances, the situation is quite different when $V$ is
(or contains) a plane or a sphere, in which case the bound goes up to $\Omega(n/\log n)$ \cite{GK2,Tao}
(see also Sheffer's survey~\cite{She} for details).

\paragraph{Incidence graph decomposition (for arbitrary points and surfaces).}
Our final main result on point-surface incidences deals with the
general setup involving a set $S$ of constant-degree algebraic
surfaces and an arbitrary set of points in $\reals^3$. The analysis
in this general setup proceeds by a recursive argument, based on the
polynomial partitioning technique of Guth and Katz~\cite{GK2}, in
which Theorem~\ref{main} plays a central role\footnote{Ideally,
applying Theorem~\ref{incth} would yield a better estimate, but,
unfortunately, we cannot control the polynomial generated by the
polynomial partitioning technique.}. This result extends a recent
result in preliminary work by the authors~\cite[Theorem 1.4]{SS16a}
from spheres to general surfaces, and extends the aforementioned result of
Zahl~\cite{Za}, for general algebraic surfaces, to the case where no
constraints are imposed on $G(P,S)$.
\begin{theorem} \label{th:incgen}
Let $P$ be a set of $m$ points in $\reals^3$, and let $S$ be a set
of $n$ surfaces from some $s$-dimensional family\footnote{Here we
use the general notion of $s$-dimensionality, not confined to points
on a variety.} $\F$ of surfaces, of constant maximum degree $E$ in
$\reals^3$. Then the incidence graph $G(P,S)$ can be decomposed as
\begin{equation} \label{gpsc}
G(P,S) = G_0(P,S) \cup \bigcup_\gamma (P_\gamma\times S_\gamma) ,
\end{equation}
where the union is now over all curves $\gamma$ of intersection of
at least two of the surfaces of $S$, and, for each such $\gamma$,
$P_\gamma = P\cap \gamma$ and $S_\gamma$ is the set (of size at least two)
of surfaces in $S$ that contain $\gamma$. Moreover, we have, for any $\eps>0$,
\begin{equation} \label{eq:mainc}
J(P,S) := \sum_\gamma \big( |P_\gamma| + |S_\gamma| \big) = O\left(
m^{\frac{2s}{3s-1}}n^{\frac{3s-3}{3s-1}+\eps} + m + n \right) ,
\quad\quad \text{and}\quad\quad |G_0(P,S)| = O(m+n) ,
\end{equation}
where the constants of proportionality depend on $\eps$, $s$, $D$,
$E$, and the complexity of the family $\F$.
\end{theorem}
As already noted, this result extends Zahl's bound~\cite{Za} to the
case where no restrictions are imposed on the incidence graph (see
the remark following Theorem~\ref{main}). Zahl's bound is the same
as ours, except for the extra factor $n^\eps$ in our bound.

We also note that Theorem~\ref{th:incgen} only applies to
$s$-dimensional families $\F$, and not to families with $k$ degrees
of freedom. The main issue here is that in Theorems~\ref{main}
and~\ref{incth}, the notion of $k$ degrees of freedom (and that of
$s$-dimensionality) is applied to the intersection curves of the
surfaces from $\F$ with some constant-degree variety, whereas here
it has to hold for the surfaces themselves in the entire
three-dimensional space. So far we are lacking a good definition of this notion that will facilitate
certain steps in the proof. See a discussion of this issue following
the proof, in Section~\ref{sec:inc}.

\subsection{The main techniques}
There are three main ingredients used in
our approach. The first ingredient, already mentioned in the context
of planar point-curve incidences, is the techniques of Pach and
Sharir~\cite{PS} (given in Theorem~\ref{th:PS}), and of Sharir and
Zahl~\cite{SZ} (Theorem~\ref{incPtCu}) concerning incidences between
points and algebraic curves in the plane. The latter bound will
be used in the analysis of incidences both between points and curves,
and between points and surfaces.

The second ingredient, relevant to the proofs of
Theorems~\ref{th:mainInc3} and~\ref{th:incgen}, is the
\emph{polynomial partitioning technique} of Guth and
Katz~\cite{GK2}, and its more recent extension by Guth~\cite{Guth},
which yields a divide-and-conquer mechanism via space decomposition
by the zero set of a suitable polynomial. This will produce
subproblems that will be handled recursively, and will leave us with
the overhead of analyzing the incidence pattern involving the points
that lie on the zero set itself. The latter step will be
accomplished by a straightforward application of Theorem~\ref{main}. We assume familiarity of the reader with these results; more details will be given in the applications of this technique in the proofs of the aforementioned theorems.

The third ingredient arises in the proof of Theorems~\ref{th:mainInc3} and
\ref{th:imprInc3}, where we argue that a ``generic'' point on a
variety $V$, that is not infinitely ruled by constant-degree curves of some given family,
as in the statement of the theorems, is incident to at most a
constant number of the given curves that are fully contained in $V$.
Moreover, we can also control the number and structural properties of ``non-generic'' points.

Before formally stating, in detail, the technical properties that we need, we
review a few notations.

Fix a constructible set $\C_0 \subset \C_{3,E}$ of irreducible
curves of degree at most $E$ in 3-dimensional space, and a
trivariate polynomial $f$. Following Guth and Zahl~\cite[Section
9]{GZ}, we call a point $p\in Z(f)$ a \emph{$(t, \C_0,
r)$-flecnode}, if there are at least $t$ curves $\gamma_1,\ldots,
\gamma_t \in \C_0$, such that, for each $i=1,\ldots, t$, (i)
$\gamma_i$ is incident to $p$, (ii) $p$ is a non-singular point of
$\gamma_i$, and (iii) $\gamma_i$ osculates to $Z(f)$ to order $r$ at
$p$. This is a generalization of the notion of a \emph{flecnodal
point}, due to Salmon~\cite[Chapter XVII, Section III]{salmon} (see
also \cite{GK2, SS4d} for more details). Our analysis requires the
following theorem. It is a consequence of the analysis of Guth and
Zahl~\cite[Corollary 10.2]{GZ}, which itself is a generalization of
the Cayley--Salmon theorem on surfaces ruled by lines (see, e.g.,
Guth and Katz~\cite{GK2}), and is closely related to
Theorem~\ref{th:gz} (also due to Guth and Zahl~\cite{GZ}). The
novelty in this theorem is that it addresses surfaces that are
\emph{infinitely ruled} by certain families of curves, where the
analysis in~\cite{GZ} only handles surfaces that are \emph{doubly
ruled} by such curves.
\begin {theorem} \label{th:ruledgz}
(a) For given integer parameters $c$ and $E$, there are constants
$c_1=c_1(c,E)$, $r=r(c,E)$, and $t=t(c,E)$, such that the following
holds. Let $f$ be a complex irreducible polynomial of degree $D \gg
E$, and let $\C_0 \subset \C_{3,E}$ be a constructible set of
complexity at most $c$. If there exist at least $c_1 D^2$ curves of
$\C_0$, such that each of them is contained in $Z(f)$ and contains
at least $c_1 D$ points on $Z(f)$ that are $(t,\C_0,r)$-flecnodes,
then $Z(f)$ is infinitely ruled by curves from $\C_0$.

\noindent{(b)} In particular, if $Z(f)$ is not infinitely ruled by curves from
$\C_0$ then, except for at most $c_1 D^2$ \emph{exceptional} curves,
every curve in $\C_0$ that is fully contained in $Z(f)$ is incident
to at most $c_1 D$ points that are incident to at least $t$ curves
in $\C_0$ that are also fully contained in $Z(f)$.
\end {theorem}

Note that, by making $c_1$ sufficiently large (specifically, choosing $c_1 > E$),
the assumption that each of the $c_1D^2$ curves in the premises of the theorem is
fully contained in $Z(f)$ follows (by B\'ezout's
theorem) from the fact that each of them contains at least
$c_1D$ points on $Z(f)$.
Although the theorem is a corollary of the work of Guth and Zahl in~\cite{GZ}, we review
(in the appendix) the machinery needed for its proof, and sketch a brief version
of the proof itself, for the convenience of the reader and in the interest of completeness.

\section{Proofs of Theorems~\protect{\ref{th:mainInc3}} and~\protect{\ref{th:imprInc3}} (points and curves)}
\label{sec:curves}

The proofs of both theorems are almost identical, and they differ in
only one step in the analysis. We will give a full proof of
Theorem~\ref{th:mainInc3}, and then comment on the few modifications
that are needed to establish Theorem~\ref{th:imprInc3}.

\paragraph{Proof of Theorem~\ref{th:mainInc3}.}
Since the family $\C$ has $k$ degrees of freedom with multiplicity
$\mu$, the incidence graph $G(P,\C)$, as a subgraph of $P\times\C$,
does not contain $K_{k,\mu+1}$ as a subgraph. The K{\H
o}v\'ari-S\'os-Tur\'an theorem (e.g., see~\cite[Section 4.5]{Mat02})
then implies that $I(P,\C) = O(mn^{1-1/k}+n)$, where the constant of
proportionality depends on $k$ (and $\mu$). We refer to this as the
\emph{naive bound} on $I(P,\C)$. In particular, when $m=O(n^{1/k})$,
we get $I(P,\C) = O(n)$. We may thus assume that $m\ge a'n^{1/k}$,
for some absolute constant $a'$.

The proof proceeds by double induction on $n$ and $m$, and establishes the bound
\begin{equation} \label{eqind}
I(P,\C) \le A \left(m^{\frac{k}{3k-2}}n^{\frac{3k-3}{3k-2}} +
m^{\frac{k}{2k-1}}n^{\frac{k-1}{2k-1}}q^{\frac{k-1}{2k-1}}+m+n\right) ,
\end{equation}
for a suitable constant $A$ that depends on $k$, $\mu$, $E$, and the complexity of $\C_0$.

The base case for the outer induction on $n$ is $n\le n_0$, for a
suitable sufficiently large constant $n_0$ that will be set later.
The bound (\ref{eqind}) clearly holds in this case if we choose $A \ge n_0$.

The base case for the inner induction on $m$ is $m\le a'n^{1/k}$,
in which case the naive bound implies that $I(P,\C) = O(n)$,
so (\ref{eqind}) holds with a sufficiently large choice of $A$.
Assume then that the bound (\ref{eqind}) holds for all sets $P'$,
$\C'$ with $|\C'|<n$ or with $|\C'|=n$ and $|P'|<m$, and let $P$ and
$\C$ be sets of sizes $|P|=m$, $|\C|=n$, such that $n>n_0$, and $m>a'n^{1/k}$.

It is instructive to notice that the two terms
$m^{\frac{k}{3k-2}}n^{\frac{3k-3}{3k-2}}$ and $m$ in (\ref{eqind})
compete for dominance; the former (resp., latter) dominates when
$m\le n^{3/2}$ (resp., $m\ge n^{3/2}$). One therefore has to treat
these two cases somewhat differently; see below and also in earlier
works~\cite{GK2,SS3d}.

\smallskip
\noindent {\bf Applying the polynomial partitioning technique.}
We construct a \emph{partitioning polynomial} $f$ for the set $\C$ of curves,
as in the recent variant of the polynomial partitioning technique, due
to Guth~\cite{Guth}. Specifically, we choose a degree

\begin{equation}
\label{eq:pocu} D = \begin{cases} c m^{\frac {k} {3k-2}}/n^{\frac 1
{3k-2}} , & \text{for $a'n^{1/k}\le m\le an^{3/2}$} ,
\\ cn^{1/2} , & \text{for $m > an^{3/2}$} ,
\end{cases}
\end {equation}

for suitable constants $c$, $a$, and $a'$ (whose values will be set
later), and obtain a polynomial $f$ of degree at most $D$, such that
each of the $O(D^3)$ (open) connected components of
$\reals^3\setminus Z(f)$ is crossed by at most $O(n/D^2)$ curves of
$\C$, where the former constant of proportionality is absolute, and
the latter one depends on $E$. Note that in both cases $1\le D\ll
n^{1/2}$, if $a$, $a'$, and $c$ are chosen appropriately.  Denote
the cells of the partition as $\tau_1, \ldots, \tau_u$, for $u =
O(D^3)$. For each $i=1,\ldots,u$, let $\C_i$ denote the set of
curves of $\C$ that intersect $\tau_i$, and let $P_i$ denote the set
of points that are contained in $\tau_i$. We set $m_i=|P_i|$ and
$n_i=|\C_i|$, for $i=1,\ldots,u$, put $m' = \sum_i m_i \le m$, and
notice that $n_i = O(n/D^2)$, for each $i$. An obvious property
(which is a consequence of the generalized version of B\'ezout's
theorem~\cite{Fu84}) is that every curve of $\C$ intersects at most
$ED+1$ cells of $\RR^3\setminus Z(f)$.

When $a'n^{1/k}\le m\le an^{3/2}$, within each cell $\tau_i$ of the
partition, for $i=1,\ldots, u$, we use the naive bound
$$
I(P_i,\C_i) = O(m_in_i^{1-1/k} + n_i) = O\left(m_i(n/D^2)^{1-1/k}
+ n/D^2 \right) ,
$$
and, summing over the $O(D^3)$ cells, we get a total of
$$
O\left( \frac{ mn^{1-1/k} }{D^{2(1-1/k)}} + nD \right).
$$
With the above choice of $D$, we deduce that the total number of
incidences within the cells is
$$
O\left(m^{\frac{k}{3k-2}}n^{\frac{3k-3}{3k-2}}\right).
$$
When $m>an^{3/2}$, within each cell $\tau_i$ of the partition we
have $n_i = O(n/D^{2})=O(1)$, so the number of incidences within
$\tau_i$ is at most $O(m_in_i) = O(m_i)$, for a total of $O(m)$ incidences.
Putting these two alternative bounds together, we get a total of
\begin{equation}
\label{eq:th6}
O\left(m^{\frac{k}{3k-2}}n^{\frac{3k-3}{3k-2}}+m\right)
\end{equation}
incidences within the cells.

\smallskip

\noindent {\bf Incidences within the zero set $Z(f)$.}
It remains to bound incidences with points that lie on $Z(f)$.
Set $P^*:= P \cap Z(f)$ and $m^*:=|P^*|=m-m'$. Let $\C^*$ denote
the set of curves that are fully contained in $Z(f)$, and set
$\C': = \C \setminus \C^*$,
$n^* := |\C^*|$, and $n' := |\C'| =n-n^*$. Since every curve of
$\C'$ intersects $Z(f)$ in at most $ED = O(D)$ points, we have (for
either choice of $D$)
\begin{equation} \label{eq:inc'}
I(P^*,\C') = O(nD) =
O\left(m^{\frac{k}{3k-2}}n^{\frac{3k-3}{3k-2}}+m\right) .
\end{equation}
Finally, we consider the number of incidences between points of
$P^*$ and curves of $\C^*$. Decompose $f$ into (complex)
irreducible components $f_1,\ldots,f_t$, for $t\le \Deg$, and assign
each point $p\in P^*$ (resp., curve $\gamma\in\C^*$) to the first
irreducible component $f_i$, such that $Z(f_i)$ contains $p$ (resp.,
fully contains $\gamma$; such a component always exists). The number of
``cross-incidences'', between points and curves assigned to
different components, is easily seen, arguing as above, to be $O(nD)$,
which satisfies our bound. In what follows, we recycle the symbols $m_i$
(resp., $n_i$), to denote the number of points (resp., curves) assigned to
$f_i$, and put $D_i=\deg(f_i)$, for $i=1,\ldots, t$. We clearly have
$\sum_i m_i = |P^*| = m^*$, $\sum_i n_i = |\C^*| = n^*$, and
$\sum_i D_i=\deg(f)=\Deg$.

For each $i=1,\ldots, t$, there are two cases to consider.

\smallskip

\noindent{\bf Case 1: $Z(f_i)$ is infinitely ruled by curves of $\C_0$.}
By assumption, there are at most $q$
curves of $\C$ on $Z(f_i)$, implying that $n_i \le q$. We project
the points of $P_i$ and the curves of $\C_i$ onto some generic
plane $\pi_0$. A suitable choice of $\pi_0$ guarantees that (i) no pair of intersection points or points of $P_i$ project to the same point, (ii) if $p$ is not incident to $\gamma$ then the
projections of $p$ and of $\gamma$ remain non-incident, (iii) no
pair of curves in $\C_i$ have overlapping projections, and (iv) no
curve of $\C_i$ contains any segment orthogonal to $\pi_0$. Moreover, the number of degrees of freedom does not change in the projection
(see Sharir et al.~\cite{SSS}). 
The number of incidences for the points and curves assigned to $Z(f_i)$ is equal to the number of incidences between the projected points and curves, which, by Theorem~\ref{th:PS}, is
$$
O\left(m_i^{\frac{k}{2k-1}}n_i^{\frac{2k-2}{2k-1}}+m_i+n_i\right) =
O\left(m_i^{\frac{k}{2k-1}}n_i^{\frac{k-1}{2k-1}}q^{\frac{k-1}{2k-1}}+m_i+n_i\right).
$$
Summing over $i=1,\ldots, t$, and using H\"older's inequality, we get the bound
$$
O\left(m^{\frac{k}{2k-1}}n^{\frac{k-1}{2k-1}}q^{\frac{k-1}{2k-1}}+m+n\right) ,
$$
which, by making $A$ sufficiently large, is at most
\begin{equation} \label{infrule}
\tfrac{A}{4}\left(m^{\frac{k}{2k-1}}n^{\frac{k-1}{2k-1}}q^{\frac{k-1}{2k-1}}+m+n\right) .
\end{equation}

\smallskip

\noindent{\bf Remark.}
This is the only step in the proof where
being of reduced dimension $s$, for $s$ sufficiently small, might
yield an improved bound (over the one in (\ref{infrule})); see
below, in the follow-up proof of Theorem~\ref{th:imprInc3}, for details.

\smallskip

\noindent{\bf Case 2: $Z(f_i)$ is not infinitely ruled by curves of
$\C_0$.} In this case, Theorem~\ref{th:ruledgz}(b) implies that
there exist suitable constants $c_1$, $t$ that depend on $E$ and on
the complexity of $\C_0$, such that there are at most $c_1 D_i^2$
exceptional curves, namely, curves that contain at least $c_1 D_i$
points that are incident to at least $t$ curves from $\C^*$.
Therefore, by choosing $c$ (in the definition of $D$) sufficiently
small, we can ensure that, in both cases (of small $m$ and large
$m$), $\sum_i D_i^2 \le (\sum_i D_i)^2 = D^2\ll n$. This allows us
to apply induction on the number of curves, to handle the
exceptional curves. Concretely, we have an inductive instance of the
problem involving $m_i$ points and at most $c_1 D_i^2 \ll n$ curves
of $\C$. By the induction hypothesis, the corresponding incidence
bound is at most
$$
A \left(m_i^{\frac{k}{3k-2}}(c_1 D_i^2)^{\frac{3k-3}{3k-2}} +
m_i^{\frac{k}{2k-1}}(c_1 D_i^2)^{\frac{k-1}{2k-1}}q^{\frac{k-1}{2k-1}}
+ m_i + c_1 D_i^2 \right) .
$$
We now sum over $i$. For the first and fourth terms, we bound each $m_i$ by $m$,
and use the fact that $\sum_i D_i^\alpha \le D^\alpha$ for any $\alpha \ge 1$.
For the second terms, we use H\"older's inequality. Overall, we get the incidence bound
$$
A \left(c_1^{\frac{3k-3}{3k-2}} m^{\frac{k}{3k-2}}(D^2)^{\frac{3k-3}{3k-2}}
+ c_1^{\frac{k-1}{2k-1}} m^{\frac{k}{2k-1}}(D^2)^{\frac{k-1}{2k-1}}q^{\frac{k-1}{2k-1}}
+ m + c_1 D^2 \right) ,
$$
which, with a proper choice of $c$ (in (\ref{eq:pocu})), can be
upper bounded by
\begin{equation} \label{eq:induct}
\tfrac{A}{4}\left(m^{\frac{k}{3k-2}} n^{\frac{3k-3}{3k-2}} +
m^{\frac{k}{2k-1}}n^{\frac{k-1}{2k-1}}q^{\frac{k-1}{2k-1}}+m+n \right) .
\end{equation}
Except for these incidences, for each $f_i$, each non-exceptional
curve in $\C^*$ that is assigned to $Z(f_i)$ is incident to at most
$c_1 D_i$ points that are incident to at least $t$ curves from
$\C^*$; the total number of incidences of this kind involving the
$n_i$ curves assigned to $Z(f_i)$ and their incident points is
$O(n_iD_i)$. Other incidences involving the non-exceptional curves
in $\C$ that are assigned to $Z(f_i)$ only involve points assigned
to $Z(f_i)$ that are incident to at most $t=O(1)$ curves from
$\C^*$; the number of such point-curve incidences is thus $O(m_i t)
= O(m_i)$. Therefore, when $Z(f_i)$ is not infinitely ruled by
curves of $\C^*$, the number of incidences assigned to $Z(f_i)$ is
$O(m_i+n_iD_i)$, plus terms that are accounted for by the induction.
Summing over these components $Z(f_i)$, we get the bound $O(m+nD)$,
plus the inductive bounds in (\ref{eq:induct}), and, choosing $A$ to
be sufficiently large, these bounds will collectively be at most
\begin{equation}
\tfrac{A}{2} \left(m^{\frac{k}{3k-2}}n^{\frac{3k-3}{3k-2}} +
m^{\frac{k}{2k-1}}n^{\frac{k-1}{2k-1}}q^{\frac{k-1}{2k-1}}+m+n\right) .
\end{equation}
In summary, by choosing $A$ sufficiently large, the number of incidences is well
within the bound of~(\ref{eqind}), thus establishing the induction step,
and thereby completing the proof.
\proofend

\paragraph{Proof of Theorem~\ref{th:imprInc3}.}
The proof proceeds by the same double induction on $n$ and $m$, and
establishes the bound, for any prespecified $\eps>0$,
\begin{equation} \label{eqind-impr}
I(P,\C) \le A m^{\frac{k}{3k-2}}n^{\frac{3k-3}{3k-2}} + A_\eps \left(
m^{\frac{2s}{5s-4}}n^{\frac{3s-4}{5s-4}}q^{\frac{2s-2}{5s-4}+\eps}+m^{2/3}n^{1/3}q^{1/3}+m+n\right) ,
\end{equation}
for a suitable constant $A$ that depends on $k$, $\mu$, $s$, $E$,
and the complexity of $\C_0$, and another constant $A_\eps$ that also depends on $\eps$. The flow of the proof is very similar
to that of the preceding proof. The main difference is in the case
where some component $Z(f_i)$ of $Z(f)$ is infinitely ruled by
curves from $\C_0$. Again, in this case it contains at most $q$ curves of $\C^*$.

We take the points of $P^*$ and the curves of $\C^*$ that are
assigned to $Z(f_i)$, and project them onto some generic plane
$\pi_0$ (the same plane can be used for all such components), as in
the proof of Theorem~\ref{th:mainInc3} and get the same properties (i)--(iv)
of the projected points and curves. Let $P_i$ and $\C_i$ denote, respectively, the set of projected points and the set of projected curves; the latter is a set of $n_i$ plane
irreducible algebraic curves of constant maximum degree\footnote{%
  A projection preserves irreducibility and does not increase the degree;
  see, e.g., Harris~\cite{Harris} for a reference to these facts.}
$DE$. Moreover, as in the preceding proof, the contribution of $Z(f_i)$ to $I(P^*,\C^*)$ is equal to the number $I(P_i,\C_i)$ of incidences between $P_i$ and $\C_i$. We can now apply
Theorem~\ref{incPtCu} to $P_i$ and $\C_i$. To do so, we first note:
\begin {lemma} \label{sames-curves}
$\C_i$ is contained in an $s$-dimensional family of curves.
\end {lemma}
\noindent{\bf Proof.}
Here it is more convenient to work over the complex field $\cplx$ (see the general remark
in the introduction).
Let $\Pi_0$ denote the projection of $\cplx^3$ 
onto $\pi_0$. Let $\C_0(f_i)$ denote the family of the curves of
$\C_0$ that are contained in $Z(f_i)$, and let $\tilde{\C}_0(f_i)$
denote the family of their projections onto $\pi_0$ (under $\pi_0$).
Define the mapping $\psi: \C_0(f_i) \to \tilde{\C}_0(f_i)$, by
$\psi(\gamma)=\Pi_0(\gamma)$, for $\gamma\in \C_0(f_i)$. By Green
and Morrison~\cite{GM} (see also~\cite[Lecture 21]{Harris} and
Ellenberg et al.~\cite[Section 2]{ESZ}), $\C_0(f_i)$ and
$\tilde{\C}_0(f_i)$ are algebraic varieties and $\psi$ is a
(surjective) morphism from $\C_0(f_i)$ to $\tilde{\C}_0(f_i)$. In
general, if $\psi: X \mapsto Y$ is a surjective morphism of
algebraic varieties, then the dimension of $X$ is at least as large
as the dimension of $Y$. Indeed, Definition 11.1~in~\cite{Harris}
defines the dimension via such a morphism, provided that it is
finite-to-one. A complete proof of the general case is given in
\cite{stack}. Therefore, $\tilde{\C}_0(f_i)$ is of dimension at most
$\dim(\C_0(f_i))=s$, and the proof of the lemma is complete. $\Box$

Applying Theorem~\ref{incPtCu} to the projected points and curves, we conclude
that the number of incidences for the points and curves assigned to $Z(f_i)$ is at most
$$
B_\eps \left(m_i^{\frac{2s}{5s-4}} n_i^{\frac{5s-6}{5s-4}+\eps} +
m_i^{2/3}n_i^{2/3} + m_i + n_i \right) \le B_\eps \left(m_i^{\frac{2s}{5s-4}}
n_i^{\frac{3s-4}{5s-4}}q^{\frac{2s-2}{5s-4}+\eps} +
m_i^{2/3}n_i^{1/3}q^{1/3} + m_i + n_i \right) ,
$$
with a suitable constant of proportionality $B_\eps$ that depends on $s$ and on $\eps$. By H\"older's
inequality, summing this bound over all such components $Z(f_i)$, we
get the bound
$$
B'_\eps \left(m^{\frac{2s}{5s-4}}n^{\frac{3s-4}{5s-4}}q^{\frac{2s-2}{5s-4}+\eps}
+ m^{2/3}n^{1/3}q^{1/3} + m + n \right) ,
$$
for another constant $B'_\eps$ proportional to $B_\eps$. By making $A_\eps$ sufficiently large, this bound is at most
$$
\tfrac{A_\eps}{4}\left(m^{\frac{2s}{5s-4}}n^{\frac{3s-4}{5s-4}}q^{\frac{2s-2}{5s-4}+\eps}
+ m^{2/3}n^{1/3}q^{1/3} + m + n \right) .
$$
The rest of the proof proceeds as the previous proof, more or less
verbatim, except that we need a more careful (albeit
straightforward) separate handling of the leading term, multiplied
by $A$, and the other terms, multiplied by $A_\eps$. The induction
step then establishes the bound in (\ref{eqind-impr}) in much the
same way as above. \proofend

\smallskip

\noindent{\bf Remarks. (1)}
As already mentioned in the introduction, the ``lower-order'' terms
$$
O\left(m^{\frac{2s}{5s-4}}n^{\frac{3s-4}{5s-4}}q^{\frac{2s-2}{5s-4}+\eps}
+ m^{2/3}n^{1/3}q^{1/3} + m + n\right)
$$
in the bound are ``best possible'' in the following sense. If the
bound in Theorem~\ref{incPtCu} were optimal, or nearly optimal, in
the worst case, for points and curves of $\C_0$ that lie in a
constant-degree surface $V$ that is infinitely ruled by such curves,
the same would also hold for the lower-order terms in the bound in
Theorem~\ref{th:imprInc3}.\footnote{Theorem~\ref{incPtCu} is
formulated, and proved in~\cite{SZ}, only for \emph{plane} curves.
Nevertheless, it also holds for curves contained in a variety $V$ of
constant degree, simply by projecting the points and curves onto
some generic plane, as done in the proofs.} This is shown by a
simple packing argument, in which we take $n/q$ generic copies of
$V$, and place on each of them $mq/n$ points and $q$ curves, so as
to obtain
$$
\Omega\left((mq/n)^{\frac{2s}{5s-4}}q^{\frac{5s-6}{5s-4}} +
(mq/n)^{2/3}q^{2/3} +mq/n + q\right)
$$
incidences on each copy, for a total of
$$
(n/q)\cdot\Omega\left((mq/n)^{\frac{2s}{5s-4}}q^{\frac{5s-6}{5s-4}}
+ (mq/n)^{2/3}q^{2/3} + mq/n + q\right) =
\Omega\left(m^{\frac{2s}{5s-4}}n^{\frac{3s-4}{5s-4}}q^{\frac{2s-2}{5s-4}}
+ m^{2/3}n^{1/3}q^{1/3} +m + n\right)
$$
incidences. (This construction works when $m>n/q$. Otherwise, the
bound is linear, and clearly best possible. Also, we assume that the
lower bound does not involve the factor $q^\eps$, to simplify the
reasoning.) In particular, this remark applies to the case of points
and circles, as discussed in Theorem~\ref{ptcirc}.

\smallskip

\noindent{\bf (2)} There is an
additional step in the proof in which the fact that $\C_0$ is of
some constant (not necessarily reduced) dimension $s'$ could lead to
an improved bound. This is the base case $m = O(n^{1/k})$, where we
use the K{\H o}v\'ari-S\'os-Tur\'an theorem to obtain a linear bound
on $I(P,\C)$. Instead, we can use the result of Fox et
al.~\cite[Corollary 2.3]{FPSSZ}, and the fact that the incidence
graph does not contain $K_{k,\mu+1}$ as a subgraph, to show that,
when $m= O(n^{1/s'})$, the number of incidences is linear. The
problem is that here we need to use the dimension $s'$ of the entire
$\C_0$, rather than the reduced dimension $s$ (which, as we recall,
applies only to subsets of $\C_0$ on a variety that is infinitely
ruled by curves of $\C_0$). Typically, as already noted, $s$ is
larger than $k$ (generally twice as large as $k$), making this
bootstrapping bound inferior to what we have. Still, in cases where
$s'$ happens to be smaller than $k$, this would lead to a further
improved incidence bounds, in which the leading term is also
smaller.

\paragraph{Rich points.}
Theorems \ref{th:mainInc3} and \ref{th:imprInc3} can easily be restated as bounding the
number of $r$-\emph{rich points} for a set $\C$ of curves with $k$ degrees of freedom
(and or reduced dimension $s$) in $\reals^3$, when $r$ is at least some sufficiently
large constant. The case $r=2$ is treated in Guth and Zahl~\cite{GZ}, and
the same bound that they obtain holds for larger values of $r$ (albeit
without an explicit dependence on $r$), smaller than the threshold in the following corollary.
\begin{corollary} \label{co:rich}
(a) Let $\C$ be a set of $n$ irreducible algebraic curves, taken from some
constructible family $\C_0$ of irreducible curves of degree at
most $E$ and with $k$ degrees of freedom (with some multiplicity
$\mu$) in $\RR^3$, and assume that no surface that is infinitely
ruled by curves of $\C_0$, or, alternatively, by curves of degree at
most $E$, contains more than $q$ curves of $\C$ (e.g., make this
assumption for all surfaces of degree at most $100E^2$). Then there
exists some constant $r_0$, depending on $k$ (and $\mu$) and on
$\C_0$, or, more generally, on $E$, such that, for any $r\ge r_0$,
the number of points that are incident to at least $r$ curves of
$\C$ (so-called \emph{$r$-rich points}) is
$$
O\left( \frac{n^{3/2}}{r^{\frac{3k-2}{2k-2}}} +
\frac{nq}{r^{\frac{2k-1}{k-1}}} + \frac{n}{r} \right) ,
$$
where the constant of proportionality depends on $k$ and $E$ (and on $\mu$).

\smallskip

\noindent
(b) If $\C_0$ is also of reduced dimension $s$, the bound on the number of $r$-rich points becomes
$$
O\left( \frac{n^{3/2}}{r^{\frac{3k-2}{2k-2}}} +
\frac{nq^{\frac{2s-2}{3s-4}+\eps}}{r^{\frac{5s-4}{3s-4}}} + \frac{n}{r} \right) ,
$$
where the constant of proportionality now also depends on $s$ and
$\eps$. (Actually, the first term comes with a constant that is
independent of $\eps$.)
\end{corollary}
\noindent{\bf Proof.}
Denoting by $m_r$ the number of $r$-rich points, the corollary is obtained by combining the
upper bound in Theorem~\ref{th:mainInc3} or Theorem~\ref{th:imprInc3} with the lower bound $rm_r$.
\proofend

The bound in (b) is an improvement, for $s=k$, when $q>r^{k+\eps'}$, for another
arbitrarily small parameter $\eps'$, which is linear in the prespecified $\eps$.
(To be more precise, this is an improvement at all only when the second term dominates the bound.)

It would be interesting to close the gap, by obtaining an $r$-dependent bound also for
values of $r$ between $3$ and $r_0$.  It does not seem that the technique in
Guth and Zahl~\cite{GZ} extends to this setup.

\section{Incidences between points and circles and similar triangles in $\reals^3$}
\label{se:sim}
We first briefly discuss the fairly straightforward proof of Theorem~\ref{ptcirc}. As already discussed in the introduction, we have $k=s=3$, for the case of circles, so we can apply Theorem~\ref{th:imprInc3} in the context of circles, and obtain the bound
$$
I(P,\C) = O\left( m^{3/7}n^{6/7} + m^{2/3}n^{1/3}q^{1/3} +
m^{6/11}n^{5/11}q^{4/11+\eps} + m + n \right) ,
$$
for any $\eps>0$, where $q$ is the maximum number of the given
circles that are coplanar or cospherical. In fact, the extension of the
planar bound (\ref{eq:impbo}) to higher dimensions, due to Aronov et
al.~\cite{AKS}, asserts that, for any set $\C$ of circles in any
dimension, we have
\begin{equation} \label{ptci}
I(P,\C) = O\left(m^{2/3}n^{2/3} + m^{6/11}n^{9/11}\log^{2/11}(m^3/n) + m + n \right) ,
\end{equation}
which is slightly better than the general bound of Sharir and Zahl~\cite{SZ}
(given in Theorem~\ref{incPtCu}). If we use this bound, instead of that in
Theorem~\ref{incPtCu}, in the proof of Theorem~\ref{th:imprInc3} (specialized for the case of circles),
we get the slight improvement (in which the two constants of proportionality are now absolute)
$$I(P,\C) = O\left( m^{3/7}n^{6/7} + m^{2/3}n^{1/3}q^{1/3} +
m^{6/11}n^{5/11}q^{4/11}\log^{2/11}(m^3/q) + m + n \right),$$ which establishes Theorem~\ref{ptcirc}.

\paragraph{The number of similar triangles.}
Theorem~\ref{ptcirc} has the following interesting application. Let $P$ be a set of $n$ points in $\reals^3$, and let $\Delta=abc$ be a fixed
given triangle. The goal is to bound the number, denoted as $S_\Delta(P)$, of triangles
spanned by $P$ and similar to $\Delta$. The best known upper bound for $S_\Delta(P)$,
obtained by Agarwal et al.~\cite{AAPS}, is $O(n^{13/6})$, and the proof that establishes this bound in \cite{AAPS} is fairly
involved. Using Theorem~\ref{ptcirc}, we obtain the following simple and fairly straightforward
improvement.
\begin{theorem} \label{thm:simi}
$S_\Delta(P) = O(n^{15/7})$.
\end{theorem}
\noindent{\bf Proof.}
Following a standard strategy, fix a pair $p,q$ of points in $P$, and consider the locus
$\gamma_{pq}$ of all points $r$ such that the triangle $pqr$ is similar to $\Delta$
(when $p,q,r$ are mapped to $a,b,c$, respectively).
Clearly, $\gamma_{pq}$ is a circle whose axis (line passing through the center of
$\gamma_{pq}$ and perpendicular to its supporting plane) passes through $p$ and $q$.
Moreover, there exist at most two (ordered) pairs $p,q$ and $p',q'$ for which
$\gamma_{pq} = \gamma_{p'q'}$. Let $\C$ denote the set of all these circles (counted
without multiplicity). Then $S_\Delta(P)$ is at most two thirds of the number $I(P,\C)$
of incidences between the $n$ points of $P$ and the $N=O(n^2)$ circles of $\C$.

By Theorem~\ref{ptcirc} we thus have
$$
S_\Delta(P) = O\left( n^{3/7}(n^2)^{6/7} + n^{2/3}(n^2)^{1/3}q^{1/3} +
n^{6/11}(n^2)^{5/11}q^{4/11}\log^{2/11} n + n^2 \right) ,
$$
where $q$ is the maximum number of circles in $\C$ that are either coplanar or cospherical.
That is, we have
\begin{equation} \label{sdp}
S_\Delta(P) = O\left( n^{15/7} + n^{4/3}q^{1/3} + n^{16/11}q^{4/11}\log^{2/11} n + n^2 \right) .
\end{equation}
We claim that $q=O(n)$. This is easy for coplanarity, because, for any fixed
plane $\pi$, each point $p\in P$ can generate at most one circle $\gamma_{pq}$
in $\C$ that is contained in $\pi$. Indeed, the axis of such a circle is perpendicular
to $\pi$ and passes through $p$. This fixes the center of $\gamma_{pq}$, and it is
easily checked that the radius is also fixed. A similar argument holds for cospherical
circles. Here too, for a fixed sphere $\sigma$, each point $p\in P$ that is not the
center $o$ of $\sigma$ can generate at most one circle $\gamma_{pq}$ in $\C$ that is
contained in $\sigma$. This is because the axis of such a circle must pass through $o$,
which fixes the center of the circle, and the radius is also fixed, as an easy
calculation shows. For $p=o$ there are at most $n-1$ additional such circles.

Hence, plugging $q=O(n)$ into (\ref{sdp}), we get $S_\Delta(P) = O(n^{15/7})$, as asserted.
$\Box$

\section{Proof of Theorem~\ref{main} (points on a variety and surfaces)} \label{sec:main}

Let $P$, $V$, $S$, $\F$, $m$, and $n$ be as in the statement of the
theorem. We first restrict the analysis to the case where $V$ is
irreducible. This involves no loss of generality, because, when $V$
is reducible, we can decompose it into its irreducible components,
assign each point of $P$ to each component that contains it, and
assign the surfaces of $S$ to all the components. This decomposes
the problem into at most $D$ subproblems, each involving an
irreducible surface, and it thus follows that the original vertex
set count is at most $D=O(1)$ times the bound for the irreducible
case. In the remainder of this section we thus assume that $V$ is
irreducible. To obtain the bound in (\ref{eq:maink}) or in
(\ref{eq:mains}) on $\sum_\gamma |P_\gamma|$, we reduce this problem
to the case of incidences between points and algebraic curves in the
plane, and then apply either Theorem~\ref{th:PS} or
Theorem~\ref{incPtCu}, as appropriate.

\paragraph{Surfaces with $k$ degrees of freedom.}
%
%
Recall that a family $\F$ of surfaces is said to have
$k$ degrees of freedom with respect to a constant-degree variety $V$, if the
family of the irreducible components of the intersection curves
$\{\sigma\cap V \mid \sigma\in\F\}$, counted without multiplicity,
has $k$ degrees of freedom, with some constant multiplicity $\mu$, as defined for curves in $\reals^3$
in Section~\ref{sec:res}.

Note that this definition means that, for any $k$ points on $V$
there are at most $\mu$ curves of the form $\sigma\cap V$, for
$\sigma\in\F$, that pass through all the points; the number of
surfaces that pass through all the points could be much larger, even
infinite. For example, spheres in $\reals^3$ have four degrees of
freedom with respect to any variety that is neither a sphere nor a
plane, because four non-cocircular points determine a unique sphere
that passes through all four, whereas four cocircular points
(over-)determine a unique circle that passes through all of them,
but an infinity of spheres with this property. Interestingly, when
$V$ is a sphere or a plane, the number of degrees of freedoms goes
down to three.

\paragraph{The reduction.}
Consider the intersection curves $\gamma_\sigma := \sigma\cap V$,
for $\sigma\in S$. These are algebraic curves of degree
$O(DE)=O(1)$. These curves are not necessarily distinct, and pairs
of distinct curves can have common irreducible components. We denote
by $\Gamma$ the multiset of the irreducible components of these
curves, where each component appears with multiplicity equal to the
number of surfaces that contain it; we also denote by $\Gamma_0$ the
underlying set of distinct irreducible components of the curves,
obtained by removing duplications from $\Gamma$.

We may assume that $V$ does not fully contain any surface of $S$.
Indeed, since $V$ is irreducible, it can contain (that is, coincide
with) at most one such surface, which contributes at most $m$ to
$\sum_\gamma |P_\gamma|$. Ignoring this surface, we have that each
$\gamma_\sigma$ is at most one-dimensional; it can be empty, and it
may have isolated points. To treat these points, we note that
each such curve has only $O(1)$ such points,\footnote{%
  The number of isolated points on a curve is easily seen to be quadratic in its degree.
  In our case, this degree is $O(DE)=O(1)$, and the claim follows.}
so the isolated points contribute a total of at most $O(n)$
incidences with their corresponding surfaces. For uniformity, we
simply record these incidences, as well as those involving the
surface fully contained in (coinciding with) $V$, if any, as trivial
(complete) bipartite graphs, with total vertex set size $O(m+n)$ on
the $P$-side, and $O(n)$ on the $S$-side.

We represent (the remainder of) $G(P,S)$ simply as the union
$\bigcup_{\gamma\in\Gamma_0} \left(P_\gamma\times S_\gamma\right)$, where,
for each $\gamma\in\Gamma_0$,
$P_\gamma = P\cap\gamma$ and $S_\gamma$ is the set of all surfaces
in $S$ that contain $\gamma$. This representation is not necessarily
edge disjoint, but a pair $(p,\sigma)$ can appear in this union at
most $O(DE)$ times, because $\sigma\cap V$ can have at most $O(DE)$
irreducible components, and $(p,\sigma)$ appears in the union once
for each of these components that contains $p$. The argument just
offered also shows that $\sum_\gamma |S_\gamma| = O(n)$. The
corresponding sum $\sum_\gamma |P_\gamma|$ (excluding the special
cases treated above, which only add $O(m+n)$ to the count) is the
number of incidences $I(P,\Gamma_0)$ between the points of $P$ and
the curves of $\Gamma_0$ (counted without multiplicity).
We therefore proceed to estimate $I(P,\Gamma_0)$.

We follow an argument very similar to the one in the proofs of
Theorems~\ref{th:mainInc3} and~\ref{th:imprInc3}; due to certain
differences, some of which are rather nontrivial, we spell it out
for clarity. We take a generic plane $\pi_0$ and project the points
of $P$ and the curves of $\Gamma_0$ onto $\pi_0$. As before, a
suitable choice of $\pi_0$ guarantees that (i) no pair of
intersection points or points of $P$ project to the same point, (ii)
if $p$ is not incident to $\gamma$ then the projections of $p$ and
of $\gamma$ remain non-incident, (iii) no pair of curves in
$\Gamma_0$ have overlapping projections, and (iv) no curve of
$\Gamma_0$ contains any segment orthogonal to $\pi_0$. Let $P^*$ and
$\Gamma_0^*$ denote, respectively, the set of projected points and
the set of projected curves; the latter is a set of $n$ plane
irreducible algebraic curves of constant maximum degree $O(DE)$ (see
a previous footnote). Moreover, $I(P,\Gamma_0)$ is equal to the
number $I(P^*,\Gamma_0^*)$ of incidences between $P^*$ and
$\Gamma_0^*$.

We now bifurcate according to which property $\F$ is assumed to
satisfy. Consider first the case where $\F$ is of reduced dimension
$s$ with respect to $V$. We have the following lemma, which is an
extension of the simpler variant given in Lemma~\ref{sames-curves}
above.
\begin {lemma} \label{sames}
$\Gamma_0^*$ is contained in an $s$-dimensional family of curves.
\end {lemma}
\noindent{\bf Proof.}
As in the preceding proof, we work over the complex field $\cplx$.
Let $\Pi_0$ denote the projection of $\cplx^3$ 
onto $\pi_0$. Define the family of curves $\Gamma_\F = \{\sigma \cap
V \mid \sigma \in \F\}$, and let $\Gamma_\F^*$ denote the family of
the projections of the curves in $\Gamma_\F$ under $\Pi_0$. Define
mappings $\varphi: \F \to \Gamma_\F$ and $\psi: \Gamma_\F \to
\Gamma_\F^*$, by $\varphi(\sigma)=\sigma \cap V$, and
$\psi(\gamma)=\Pi_0(\gamma)$. As above, by~\cite{ESZ,GM, Harris},
$\Gamma_\F$ and $\Gamma_\F^*$ are algebraic varieties and $\psi$ is
a morphism from $\Gamma_\F$ to $\Gamma_f^*$. By~Fulton~\cite[Section
3.4]{Fu84}, $\varphi$ is a morphism from $\F$ to $\Gamma_\F$, and
since both $\varphi$ and $\psi$ are surjective, it follows that
their composition $\Phi = \psi \circ \varphi$ is a surjective
morphism from $\F$ to $\Gamma_\F^*$. As in the proof of
Lemma~\ref{sames-curves}, the definition in Harris~\cite[Definition
11.1]{Harris}, and its extension in~\cite{stack}, imply that the
dimension of $\F$ is at least as large as the dimension of
$\Gamma_\F^*$. Therefore, $\Gamma_\F^*$ is of dimension at most
$\dim(\F)=s$, and the proof of the lemma is complete. $\Box$

\smallskip

\noindent{\bf Remark.} It might be the case that $\Gamma_0^*$ is of
smaller dimension than $s$. As will follow from the proof, the
incidence bound depends on the dimension of $\Gamma_0^*$ and not on
the dimension of $\F$, so the bound will improve if $\Gamma_0^*$ is
indeed of smaller dimension.

In addition, the curves of $\Gamma_0^*$ are of constant maximum
degree. Applying Theorem~\ref{incPtCu} to our setup, we get the
bound in (\ref{eq:mains}). Adding the counts obtained separately for
the preceding special cases completes the proof of
Theorem~\ref{main} when $\F$ is of reduced dimension $s$ with
respect to $V$.

Consider next the case where $\F$ has $k$ degrees of freedom with respect to $V$. Then
Theorem~\ref{th:PS} is applicable to $P^*$ and $\Gamma_0^*$, and we have
$$
I(P^*,\Gamma_0^*) =
O\left(m^{\frac{k}{2k-1}}n^{\frac{2k-2}{2k-1}}+m+n\right) ,
$$
and adding to this the bounds obtained in the other cases yields the bound asserted in (\ref{eq:maink}).
$\Box$

\section{Proof of Theorem~\ref{incth} (points on a variety and general surfaces)} \label{sec:incth}

Let $P$, $V$, $S$, $\F$, $k$, $\mu$, $s$, $m$, and $n$ be as in the
statement of the theorem. We first restrict the analysis to the case
where $V$ is irreducible. The general case can be handled, as in the
case of Theorem~\ref{main}, by repeating the analysis to each of the
$O(1)$ irreducible components of $V$, and summing up the resulting
bounds, to obtain the same asymptotic bound (multiplied by an extra
factor of $D=O(1)$).

Let then $f$ be an irreducible complex polynomial such that $V=Z(f)$,
and let $\C$ denote the set of irreducible curves that (i) are fully
contained in $V$, (ii) are contained in at least two surfaces of
$S$, and (iii) contain at least one point of $P$.
By B\'ezout's theorem \cite{Fu84} and condition (ii), we have
$\deg(\gamma)\le E^2$ for each $\gamma\in \C$. For each curve
$\gamma\in \C$ we form the bipartite subgraph $P_\gamma \times
S_\gamma$ of $G(P,S)$, where $P_\gamma = P\cap \gamma$ (the actual
sets $P_\gamma$ for some of the curves will be smaller---see
below), and $S_\gamma$ is the set of the surfaces of $S$ that
contain $\gamma$. To estimate $\sum_\gamma |P_\gamma|$ and
$\sum_\gamma |S_\gamma|$, we argue as follows. First, $\sum_\gamma
|P_\gamma|$ is the number of incidences between the points of $P$
and the curves of $\C$, counted without multiplicity. By assumption,
$V$ is not infinitely ruled by the irreducible components of the
intersection curves of pairs of surfaces from $\F$, and $\C$ is
contained in this family of curves. To apply
Theorem~\ref{th:ruledgz}, it remains to argue that $\C$ is
constructible. Although the proof of this property is not too hard,
it is rather technical, and we will present it in
Lemma~\ref{c-const} in the Appendix. We then conclude that, for a
suitable constant $t=t(E, c)$ (where $c$ is the complexity of the
family $\C$), except for possibly $O(D^2)$ exceptional curves, every
curve in $\C$ contains only $O(D)$ points that are incident to at
least $t$ curves of $\C$.

Consider first incidences between the non-exceptional curves in $\C$
and the ``rich'' points (those incident to at least $t$ curves of
$\C$). Each surface $\sigma \in S$ intersect $V$ in a curve of
degree at most $DE$, and can therefore contain at most $O(DE)$
curves from $\C$. Each of these curves contains at most $O(D)$
$t$-rich points, for a total of $O(D^2E)$ incidences for each
surface $\sigma \in S$, so the overall number of incidences of this
kind is $O(nD^2E)=O(n)$. Note that this is a bound on the actual
number of point-surface incidences. We include the incident pairs of
this kind in $G_0(P,S)$, and the resulting bound $O(n)$ is clearly
subsumed by the asserted bound in (\ref{eq:main2k}) or
(\ref{eq:main2s}).

Removing these edges from the complete bipartite decomposition, the
remaining incidences counted in $I(P,\C)$ are estimated as follows.
The number of incidences with the $t$-poor points (each lying on at
most $t$ curves of $\C$) is at most $mt=O(m)$, and the number of
incidences between the $O(D^2)$ exceptional curves and the points of
$P$ is at most $O(mD^2)=O(m)$, for a total of $O(m)$ incidences.

In summary, the complete bipartite decomposition that we end up with
is of the form $\bigcup_\gamma P'_\gamma \times S_\gamma$, where
$\gamma$ ranges over the curves of $\C$, and (i) for each of the
$O(D^2)$ exceptional curves $\gamma$ we have $P'_\gamma = P_\gamma$,
and (ii) for each of the non-exceptional curves $\gamma$,
$P'_\gamma$ is the set of the $t$-poor points of $P$ that lie on
$\gamma$. We thus obtain that $\sum_{\gamma\in \C}
|P'_\gamma|=O(m)$. As $V$ is irreducible and does not contain any
of the surfaces $\sigma\in S$ (except possibly for at most one,
which then coincides with $V$ and which we may ignore, as before),
the preceding argument implies again that each $\sigma \in S$
generates at most $DE$ curves of $\C$, so $\sum_{\gamma} |S_\gamma|
\le nDE = O(n)$. This gives us the complete bipartite graph
decomposition portion of the representation in (\ref{gps1}), which
satisfies all the properties asserted in the theorem.


Let $G_0(P,S)$ denote the remainder of the incidence graph. For the
moment, ignore the pairs involving the $t$-rich points on the curves of
$\C$, which are also part of the final $G_0(P,S)$. For each
$\sigma\in S$, put $\gamma_\sigma := (\sigma\cap V) \setminus
\bigcup \C$. As just noted, each $\gamma_\sigma$ is at most
one-dimensional (i.e., a curve). By construction, it does not
contain any curve in $\C$, and it might also be empty (for this or
for other reasons). Note that if $\sigma\cap V$ does contain a curve
$\gamma$ in $\C$, then the incidences between $\sigma$ and the
points of $P$ on $\gamma$ are all already recorded in
$P'_\gamma\times S_\gamma$, or are the $O(n)$ special incidences with
$t$-rich points on the curves in $\C$, so ignoring them is ``safe''.
Finally, we may ignore the isolated points of $\gamma_\sigma$,
because, as already argued, each curve $\gamma_\sigma$ can contain
at most $O(1)$ such points, which contribute a total of at most
$O(n)$ incidences with their corresponding surfaces. Let $G_0(P,S)$
continue to denote the remaining portion of $G(P,S)$, after pruning
away all the incidences already accounted for. Put $I_0(P,S) =
|G_0(P,S)|$.

Let $\Gamma$ denote the set of the $n$ curves $\gamma_\sigma$, for
$\sigma\in S$ (and notice that this time it is an actual set, not a
multiset). The curves of $\Gamma$ are algebraic curves of degree at
most $DE$, and, as is easily checked, any pair of curves
$\gamma_\sigma$, $\gamma_{\sigma'}\in \Gamma$ intersect in at most
$\min\{DE^2,E^4\}=O(1)$ points.

Note that $I_0(P,S)$ is equal to the number $I(P,\Gamma)$ of
incidences between the points of $P$ and the curves of $\Gamma$. To
bound the latter quantity, we proceed exactly as in the preceding
proof, bifurcating according to whether $\F$ is of reduced dimension
$s$ with respect to $V$ or has $k$ degrees of freedom with respect
to $V$. In both cases we project the points and curves onto some
generic plane $\pi_0$, and bound the number of incidences between
the projected points and curves, using either Theorem~\ref{th:PS}
(for families with $k$ degrees of freedom) or Theorem~\ref{incPtCu}
(for $s$-dimensional families), obtaining the respective bounds
asserted in (\ref{eq:main2k}) or in (\ref{eq:main2s}). $\Box$

\section{Distinct and repeated distances in three dimensions} \label{sec:dd3}

In this section we prove Theorems~\ref{dd3} and \ref{und}, the
applications of Theorems~\ref{main} and \ref{incth} to distinct and
repeated distances in three dimensions; see also our earlier
work~\cite{SS16a} that handles these problems in a somewhat
different manner. The theorems present four results, in each of
which the problem is reduced to one involving incidences between
spheres and points on a surface $V$. However, except for
Theorem~\ref{dd3}(b), the spheres that arise in the other three
cases are restricted, by requiring their centers to lie on $V$ and /
or to have a fixed radius. This makes the number of degrees of
freedom (with respect to the variety) and the dimensionality of the
corresponding families of spheres go down to $3$ or $2$. The case of
two degrees of freedom (in Theorem~\ref{und}(a)) is the simplest,
and requires very little of the machinery developed here (see
below). The cases of three degrees of freedom (in
Theorem~\ref{dd3}(a) and Theorem~\ref{und}(b)) yield improved
``in-between'' bounds.

\smallskip

\noindent{\bf Proof of Theorem~\ref{dd3} (distinct distances).} We
will first establish the more general bound in (b); handling (a)
will be done later, in a similar, somewhat simpler manner.

\smallskip

\noindent{\bf (b)}
Let $t$ denote the number of distinct distances in $P_1\times P_2$.
For each $q\in P_2$, draw $t$ spheres centered at $q$ and having as radii
the $t$ distinct distances. We get a collection $S$ of $nt$ spheres, a set
$P_1$ of $m$ points on $V$, which we relabel as $P$, to simplify the notation,
and exactly $mn$ incidences between the points of $P$ and the spheres of $S$.

Let $\C$ denote the set of intersection circles of pairs of spheres
from $S$ that are contained in $V$, counted without multiplicity; we
keep in $\C$ only circles that contain points of $P_1$. For each
$\gamma\in\C$, let $\mu(\gamma)\ge 2$ denote its multiplicity,
namely the number of spheres containing $\gamma$; note that
$\mu(\gamma)$ is equal to the number of points of $P_2$ that lie on
the \emph{axis} of $\gamma$, namely the line that passes through the
center of $\gamma$ and is orthogonal to the plane containing
$\gamma$. The maximum possible multiplicity of a circle is at most
$2t$, because the distances of the corresponding centers in $P_2$ to
the points of $P_1\cap\gamma\ne\emptyset$ are all distinct, up to a
possible multiplicity of $2$. For each $k$, let $\C_k$ (resp.,
$\C_{\ge k}$) denote the subset of circles in $\C$ of multiplicity
exactly (resp., at least) $k$, and put $N_k := |\C_k|$, $N_{\ge k}
:= |\C_{\ge k}|$.

In order to effectively apply the bound in Theorem~\ref{main:sph},
we first have to control the term $\sum_\gamma |P_{\gamma}|\cdot
|S_{\gamma}|$, which arises since we deal here with the actual
number of incidences. Specifically, we claim that most of the $mn$
incidences do not come from this part of the incidence graph, unless
$t=\Omega(n)$. Indeed, write this sum as
$$
\sum_\gamma |P_{\gamma}|\cdot |S_{\gamma}| = \sum_{k\ge 2} k\sum_{\gamma\in \C_k} |P_\gamma| .
$$
Putting $E_k := \sum_{\gamma\in \C_k} |P_\gamma|$, and
$E_{\ge k} := \sum_{\gamma\in \C_{\ge k}} |P_\gamma|$, we then have
$$
\sum_\gamma |P_{\gamma}|\cdot |S_{\gamma}| = \sum_{k\ge 2} kE_k =
2E_{\ge 2} + \sum_{k\ge 3} E_{\ge k} .
$$
By Theorem~\ref{incth}, we have $E_{\ge k} = O(m)$, so we
have $$
\sum_\gamma |P_{\gamma}|\cdot |S_{\gamma}| = 2E_{\ge 2} + \sum_{k\ge
3} E_{\ge k} = O\left( m + \sum_{k=3}^{2t} m \right) = O(mt) .
$$
If this would have accounted for more than, say, half the
incidences, we would get $mn = O(mt)$, or
$t=\Omega(n)$, as claimed, and then (a much
larger lower bound than) the bound in the theorem would follow. We
may thus ignore this term, and write
$$
mn = O\left( m^{1/2}(nt)^{7/8+\eps} + m^{2/3}(nt)^{2/3} + m + nt \right) ,
$$
for any $\eps>0$, or
$$
t = \Omega \left( \min\left\{ m^{4/(7+8\eps)}n^{(1-8\eps)/(7+8\eps)},\; m^{1/2}n^{1/2},\; m \right\} \right) ,
$$
which, by replacing $\eps$ by another, still arbitrarily small
$\eps'$, becomes the bound asserted in the theorem (and is also
smaller than the bound for the complementary situation treated
above).

\smallskip

\noindent{\bf (a)} Here we are in a more favorable situation,
because the spheres in $S$ belong to a three-dimensional family of
surfaces---a family that can be represented simply as
$V\times\reals$. We can therefore apply Theorem~\ref{incth} with
dimensionality $s=3$ (which is the actual dimensionality of the
family, not the reduced one with respect to $V$), arguing first
that, as in the proof of (b), we may ignore the term $\sum_\gamma
|P_\gamma|\cdot |S_\gamma|$ in the bound on $I(P,S)$, which is
negligible unless $t=\Omega(n)$. We thus get the inequality
$$
n^2 = O\left( n^{6/11}(nt)^{9/11+\eps} + n^{2/3}(nt)^{2/3} + nt \right) ,
$$
for any $\eps>0$, which yields $t=\Omega(n^{7/(9+11\eps)})$, which, by replacing $\eps$, as in
the proof of (b), can easily be massaged into the bound asserted in the theorem.
$\Box$

\smallskip

\noindent{\bf Proof of Theorem~\ref{und} (repeated distances).}
Consider (a) first. Following the standard approach to problems
involving repeated distances, we draw a unit sphere $s_p$ around
each point $p\in P$, and seek an upper bound on the number of
incidences between these spheres and the points of $P$; this latter
number is exactly twice the number of unit distances determined by
$P$.

This instance of the problem has several major advantages over the
general analysis in Theorem~\ref{main:sph}. First, in this case the
incidence graph $G(P,S)$ cannot contain $K_{3,3}$ as a subgraph,
eliminating altogether the complete bipartite graph decomposition in
(\ref{gps-sph}) (or, rather, bounding the overall number of edges in
these subgraphs by $O(n)$).

More importantly, the family of the fixed-radius spheres whose centers
lie on $V$ is $2$-dimensional and has two degrees of freedom, which
leads to the standard Szemer\'edi-Trotter-like bound
$I(P,S) = O(|P|^{2/3}|S|^{2/3}+|P|+|S|) = O(n^{4/3})$, so the number
of repeated distances in this case is $O(n^{4/3})$, as claimed.
(The last bound can be obtained by applying Theorem~\ref{main},
but it can also be obtained more directly, e.g., via Sz\'ekely's technique~\cite{Sze}.)

We remark that the last bound does not use much of the machinery developed
in this paper. Still, we are not aware of any previous claim of the bound
for the general case of points on a surface; see Brass et al.~\cite[Section 5.2]{BMP}
and Brass~\cite{Bra} for a discussion of closely-related problems.

\smallskip

\noindent We now consider (b). Again, we reduce the problem to that
of bounding the number of incidences between the $m$ points of
$P_1$, which lie on $V$, and the $n$ unit spheres centered at the
points of $P_2$. Here too the overall number of edges in the
complete bipartite graph decomposition is $O(m+n)$, so we can ignore
this part of the bound.

In this case, the family of unit spheres is $3$-dimensional.
Applying the same reasoning as in the proof of Theorem~\ref{main},
we conclude that the number of unit distances in this case is
$$
O\left( m^{6/11}n^{9/11+\eps} + m^{2/3}n^{2/3} + m + n \right) ,
$$
for any $\eps>0$, as claimed. $\Box$

\smallskip

%

\section{Proof of Theorem~\protect{\ref{th:incgen}} (surfaces and arbitrary points)} \label{sec:inc}

The proof establishes the bound in (\ref{eq:mainc}), via induction
on $m$, with a prespecified fixed parameter $\eps>0$. Concretely, we
claim that, for any such choice of $\eps>0$, we can write
$$
G(P,S) = G_0(P,S)\cup \bigcup_{\gamma\in \Gamma_0} (P_\gamma\times S_\gamma) ,
$$
where $\Gamma_0$ is a collection of distinct constant-degree
irreducible algebraic curves in $\reals^3$, and, for each
$\gamma\in\Gamma_0$, $P_\gamma = P\cap \gamma$ and $S_\gamma$ is the
set of surfaces in $S$ that contain $\gamma$. We only include in
$\Gamma_0$ curves $\gamma$ with $|S_{\gamma}|\ge 2$ (as the cases
where $|S_{\gamma}| = 1$ will be ``swallowed'' in $G_0(P,S)$), so
each $\gamma \in \Gamma_0$ is an irreducible component of an
intersection curve of at least two surfaces from $S$, and its degree
is therefore at most $E^2$. We then claim that
\begin{equation} \label{eq:ind}
J(P,S) := \sum_{\gamma\in\Gamma_0} \big( |P_\gamma| + |S_\gamma|
\big) \le A\left( m^{\frac {2s} {3s-1}} n^{\frac {3s-3}{3s-1}+\eps}
+ m + n \right) ,
\end{equation}
$$
\text{and} \quad \quad |G_0(P,S)| \le A(m+n),
$$
for a suitable constant $A$ that depends on $\eps$, $s$, $E$, $D$,
and the complexity of the family $\F$.

The base cases are when $m\le n^{1/s}$ or when $m\le m_0$, for some
sufficiently large constant $m_0$ that will be set later. Consider
first the case $m\le n^{1/s}$. Note that in this case the right-hand
side of~(\ref{eq:ind}) is $O(n^{1+\eps})$. We will actually
establish the bound $O(n)$ on both $J(P,S)$ and $|G_0(P,S)|$, as
follows. Since the surfaces of $S$ come from an $s$-dimensional
family, a suitable extension of the analysis in Sharir and
Zahl~\cite[Lemma 3.2]{SZ} shows that there exists an $s$-dimensional
real parametric space $\reals^s$, and a duality mapping that sends
each surface $\sigma\in S$ to a point $\sigma^*\in\reals^s$, and
sends each point $p\in P$ to a constant-degree algebraic
hypersurface $p^*$ in that space, so that if $p$ is incident to
$\sigma$ then $\sigma^*$ is incident to $p^*$. This holds with the
exception of at most $O(1)$ `bad' points in $P$ and at most $O(1)$
`bad' surfaces in $S$, and the constants depend on $s$, $E$, and the
complexity of $\F$. The contribution to $J(P,S)$, or rather to
$I(P,S)$, of the bad points and surfaces is only $O(m+n)$, so we can
ignore it, or, rather, place these incidences in the remainder
subgraph $G_0(P,S)$.

Construct the arrangement of the dual surfaces $p^*$,
for $p\in P$. Its complexity is $O(m^s) = O(n)$, and this bound also
holds if we count each face of the arrangement, of any dimension,
with multiplicity equal to the number of surfaces that contain it.
For each such (relatively open) face\footnote{%
  Technically, rather than considering individual faces $f$, we should
  consider the full varieties that contain these faces and are obtained
  by intersecting subsets of the surfaces $p^*$, where each such intersection
  might contain many faces $f$. However, in doing so, we want to exclude
  faces that lie on such an intersection and are of dimension smaller than that of the intersection (because they lie on other surfaces too), and treat them separately.
  To simplify the presentation, we ignore this modification, which does not
  affect the asymptotic bound that is derived here.}
$f$, we form the complete bipartite graph $P_f\times S_f$, where $S_f$ is the set of all
surfaces $\sigma$ such that $\sigma^*\in f$, and $P_f$ is the set of
all points $p$ whose dual surface $p^*$ contains $f$.
We have $\sum_f |S_f| \le n$, and $\sum_f |P_f| = O(m^s) = O(n)$.
Back in primal space, if $|S_f|\ge 2$, then all the points of $P_f$
lie in the intersection $\gamma_f:=\bigcap_{\sigma \in S_f} \sigma$, which, by
B\'ezout's theorem, is either one-dimensional, i.e., a curve in $\Gamma_0$
as in the theorem, or a discrete set of at most $E^3$
points. In the latter case $|P_f|\le E^3=O(1)$, implying that
$\sum_f |P_f||S_f|= O(\sum_f |S_f|)=O(n)$. Similarly, if $|S_f|=1$,
then we also have $\sum_f |P_f||S_f|= O(\sum_f |P_f|)=O(m^s)=O(n)$.
Clearly, $\sum_f |P_f||S_f|$, over faces $f$ for which either $|S_f|\le 1$ or
$\gamma_f$ is discrete, counts the number of incidences between $P$ and $S$
that fall into these special cases, so the number of these incidences is only $O(n)$.
We are left with a portion of the incidence graph that can be written as
the union of complete bipartite graphs $\bigcup_f P_{f}\times S_{f}$, over faces
$f$ for which $|S_f|\ge 2$ and $\gamma_f$ is a curve in $\Gamma_0$. This union
is of the form asserted in the theorem, and the corresponding
$J(P,S)$ is $O(n)$. The asserted bound thus holds by choosing $A$
sufficiently large.

The case $m\le m_0$ follows easily (since in this case we have
$I(P,S) \le m_0n$) if we choose $A$ sufficiently large. This holds
for any choice of $m_0$ (and a corresponding choice of $A$); the
value that we choose is specified later.

Suppose then that (\ref{eq:ind}) holds for all sets $P'$, $S'$, with
$|P'|<m$, and consider the case where the sets $P,S$ are of
respective sizes $m,n$, and we have $m>n^{1/s}$ and $m>m_0$.

Before continuing, we also dispose of the case $m\ge n^3$. In this
case we consider the arrangement $\A(S)$ (in $\reals^3$) of the
surfaces in $S$. The complexity of $\A(S)$ is $O(n^3)=O(m)$. More
precisely, this bound holds, and is asymptotically tight, for
surfaces in general position. In our case, $S$ is likely not to be
in general position, and then the complexity of $\A(S)$ might be
smaller, because vertices and edges might be incident to many
surfaces. Nevertheless, if we count each vertex and edge of $\A(S)$
with its multiplicity, we still get the complexity upper bound
$O(n^3)$. (Here we reason in complete analogy with the dual
$s$-dimensional case treated above.) This means that the number of
incidences with points that are either vertices or lie on the
(relatively open) 2-faces of $\A(S)$ is $O(n^3+m)=O(m)$. Incidences
with points that lie on the (relatively open) edges of $\A(S)$ (note
that each such edge is a portion of some curve of intersection
between at least two surfaces of $S$) are recorded, as usual, by a
complete bipartite graph decomposition $\bigcup_\gamma
(P_\gamma\times S_\gamma)$, where the curves $\gamma$ are as
stipulated in the theorem, and where, as just argued, we have
$\sum_\gamma |P_\gamma| \le m$ and $\sum_\gamma |S_\gamma| = O(n^3)
= O(m)$. This implies that (\ref{eq:ind}) holds in this case. Thus,
in what follows, we assume that $m\le n^3$. Since we also assume
that $m>m_0$, we have $n\ge m^{1/3} > n_0:= m_0^{1/3}$.
\paragraph{Applying the polynomial partitioning technique.}
We fix a sufficiently large constant parameter $D \ll m^{1/3}$,
whose concrete choice will be specified later, and apply the
polynomial partitioning technique of Guth and Katz~\cite{GK2}. We
obtain a polynomial $f\in\reals[x,y,z]$ of degree at most $D$, whose
zero set $Z(f)$ partitions 3-space into $O(D^3)$ (open) connected
components (cells), and each cell contains at most $O(m/D^3)$
points. By duplicating cells if necessary\footnote{By using Guth's
recent technique for partitioning sets of varieties~\cite{Guth},
already mentioned earlier, we could do without this cell duplication
step.}, we may also assume that each cell is crossed by at most
$O(n/D)$ surfaces of $S$; this duplication keeps the number of cells
$O(D^3)$ (because each surface crosses only $O(D^2)$ cells,
a well known property that follows, e.g., from Warren's theorem~\cite{War}).\footnote{%
  In actuality, the bound is $O(D^2E^2)$, because the intersection curve $\sigma\cap Z(f)$,
  for any $\sigma\in S$, is of degree at most $DE$. Here we treat $E$ as a much smaller
  quantity than $D$, and bear in mind that the relevant constants may depend on $E$.}
That is, we obtain at most $aD^3$ subproblems, for some absolute
constant $a$, each associated with some cell of the partition, so
that, for each $i\le aD^3$, the $i$-th subproblem involves a subset
$P_i\subset P$ and a subset $S_i\subset S$, such that $m_i:=|P_i|\le
b_0m/D^3$ and $n_i:=|S_i|\le bn/D$, for another absolute constant
$b_0$ and a constant $b$ that depends on $E$. Set $P_0:= P\cap Z(f)$
and $P'=P\setminus P_0$. We have
\begin{equation} \label{eqsplit}
J(P,S) \le J(P_0,S) + J(P',S) .
\end{equation}
We first bound $J(P_0,S)$. Similarly to what was done earlier,
decompose $Z(f)$ into its $O(D)$ irreducible components, assign each
point of $P_0$ to every component that contains it, and assign the
surfaces of $S$ to all components. We now fix a component, and bound
the vertex count in the complete bipartite graph decomposition
involving incidences between the points and surfaces assigned to
that component; $J(P_0,S)$ is at most $D$ times the bound that we
get. We may therefore assume that $Z(f)$ is irreducible. Since
$\deg(Z(f))\le D$ is a constant, Theorem~\ref{main} implies that we can write $G(P_0,S)$ as $\bigcup_\gamma P_{0 \gamma} S_\gamma$, over a suitable set of curves $\gamma \subset V$, with $P_{0 \gamma} = P_0 \cap \gamma$ and $S_\gamma$ is the set of surfaces in $S$ containing $\gamma$, and we have
$$
J(P_0,S) = O\Big(m^{\frac{2s}{5s-4}} n^{\frac{5s-6}{5s-4}+\eps} +
m^{2/3}n^{2/3} + m + n\Big) ,
$$
where the constant of proportionality depends on $\eps$, $s$, $D$,
$E$, and the complexity of $\F$. (Here, as in Theorem~\ref{main},
the contribution of $S$ to this bound, namely to $\sum_\gamma |S_\gamma |$, is only $O(n)$. Unfortunately,
this no longer holds in the estimation of $J(P',S)$, given, so we
ignore this improvement). As is easily checked, this bound is
subsumed in~(\ref{eq:ind}) for $m\ge n^{1/s}$ (and $s\ge 3$), if we
choose $A$ sufficiently large (so here, as in all the other steps,
$A$ depends on all the parameters just listed).
Finally, we estimate
$$
J(P',S) \le \sum_{i=1}^{aD^3} J(P_i,S_i) .
$$
By the induction hypothesis, we get, for each $i$,
$$
J(P_i,S_i) \le A\left(
m_i^{\frac{2s}{3s-1}}n_i^{\frac{3s-3}{3s-1}+\eps} + m_i + n_i
\right) .
$$
Summing this over $i$, we get
\begin{align*}
J(P',S) & \le A\cdot aD^3 \left( (b_0m/D^3)^{\frac{2s}{3s-1}} (bn/D)^{\frac{3s-3}{3s-1}+\eps} + (b_0m/D^3) + (bn/D) \right) \\
& =
\frac{Aab_0^{\frac{2s}{3s-1}}b^{\frac{3s-3}{3s-1}+\eps}}{D^{\eps}}
m^{\frac{2s}{3s-1}}n^{\frac{3s-3}{3s-1}+\eps} + Aab_0m + AabD^2n .
\end{align*}
We note that $m^{\frac{2s}{3s-1}}n^{\frac{3s-3}{3s-1}+\eps} \ge n^\eps\cdot m$
and $m^{\frac{2s}{3s-1}}n^{\frac{3s-3}{3s-1}+\eps} \ge n^\eps\cdot n$
for $n^{1/s}\le m\le n^3$. We choose $D$ sufficiently large so that
$D^{\eps} \ge 4Aab_0^{\frac{2s}{3s-1}}b^{\frac{3s-3}{3s-1}+\eps}$,
and then the bound is at most
$$
\left( \frac{A}{4} + \frac{Aab_0}{n^\eps} + \frac{AabD^2}{n^\eps}
\right) m^{\frac{2s}{3s-1}}n^{\frac{3s-3}{3s-1}+\eps} .
$$
Choosing $n_0$ (that is, $m_0$) sufficiently large, so that
$n_0^\eps \ge 4a\cdot\max\{b_0,bD^2\}$, we ensure that, for $n\ge n_0$,
$$
\frac{A}{4} + \frac{Aab_0}{n^\eps} + \frac{AabD^2}{n^\eps} \le
\frac{A}{4} + \frac{A}{4} + \frac{A}{4} = \frac{3A}{4} .
$$
Adding the bounds for $J(P_0,S)$, and choosing $A$ sufficiently
large, we get that (\ref{eq:ind}) holds for $P$ and $S$. This
establishes the induction step and thereby completes the proof.
$\Box$

\smallskip

\noindent{\bf Remark.} The dimensionality $s$ of $S$ is used in the
proof in two different steps, once in establishing a linear bound
when $m < n^{1/s}$, and once in deriving the bound on $J(P_0,S)$,
using Theorem~\ref{main}. As in the remark following
Theorem~\ref{main}, the values of $s$ used in these two steps need
not be the same, and the latter one is typically smaller (because it
is reduced, with respect to intersection curves with a variety).
Unfortunately, this in itself does not lead to an improvement in the
bound, because the leading term $m^{\frac {2s} {3s-1}} n^{\frac
{3s-3}{3s-1}+\eps}$ depends on the former value of $s$. If this
value of $s$ could also be improved, say by additional assumptions
on the points and/or the surfaces, the bound in the theorem would
improve too.


%

\section{\bf Discussion}

In this paper we have made significant progress on major incidence
problems involving points and curves and points and surfaces in
three dimensions. We have also obtained several applications of
these results to problems involving repeated and distinct distances
in three dimensions, with significantly improved lower and upper
bounds, in cases where the points, or in the bipartite versions, the
points in one of the two given sets, lie on a constant-degree
algebraic surface.

The study in this paper raises several interesting open problems.
\smallskip

\noindent{\bf (i)} A long-standing open problem is that of
establishing the lower bound of $\Omega(n^{2/3})$ for the number of
distinct distances determined by a set of $n$ points in $\reals^3$,
without assuming them to lie on a constant-degree surface. The best
known lower bound, close to $\Omega(n^{3/5})$, which follows from
the work of Solymosi and Vu~\cite{SoVu}, still falls short of this
bound.
In the present study we have obtained some partial results (with
better lower bounds) for cases where (all or some of) the points do
lie on such a surface. We hope that some of the ideas used in this
work could be applied in more general contexts, or in other special
situations.

\smallskip

\noindent{\bf (ii)} Another major long-standing open problem is that of
improving the upper bound $O(n^{3/2})$, established in
\cite{KMSS,Za}, on the number of unit distances determined by a set
of $n$ points in $\reals^3$, again without assuming (all or some of) them to lie on a
constant-degree surface. It would be interesting to make progress on this problem.

\smallskip

\noindent{\bf (iii)} As remarked above, a challenging open problem
is to characterize all the surfaces that are infinitely ruled by
algebraic curves of degree at most $E$ (or by certain classes
thereof), extending the known characterizations for lines and
circles. A weaker, albeit still hard problem is to reduce the upper
bound $100E^2$ on the degree of such a surface, perhaps all the way
down to $E$, or at least to $O(E)$.

\smallskip

\noindent{\bf (iv)} It would also be interesting to find additional
applications of the results of this paper, like the one with an
improved bound on the number of similar triangles in $\reals^3$,
given in Section~\ref{se:sim}. One direction to look at is the
analysis of other repeated patterns in a point set, such as
higher-dimensional congruent or similar simplices, which can
sometimes be reduced to point-sphere incidence problems; see
\cite{AAPS,AgS}.

\smallskip

\noindent{\bf (v)} Concerning Theorem~\ref{th:incgen}, we note that
it is stated only for families $\F$ of surfaces of a given
dimensionality $s$. It would be interesting to obtain a variant in
which we assume instead that $\F$ has $k$ degrees of freedom, after
one comes up with a definition of this notion that is both (a)
natural and simple to state, and (b) makes (a suitable variant of)
the analysis work. We are currently studying such an extension.

\noindent{\bf (vi)} A potentially weak issue in our analysis,
manifested in the proof of all our main theorems, is that in order
to bound the number of incidences between points and curves on some
variety $V$ of constant degree, we project the points and curves on
some generic plane and use a suitable planar bound, from
Theorem~\ref{th:PS} or Theorem~\ref{incPtCu}, to bound the number of
incidences between the projected points and curves. It would be very
interesting if one could obtain an improved bound, exploiting the
fact that the points and curves lie on a variety $V$ in $\reals^3$,
under suitable (natural) assumptions on $V$.

\smallskip

\noindent{\bf (vii)} Finally, it would be challenging to extend the
results of this paper to higher dimensions.


\appendix

\section{On surfaces ruled by curves} \label{ap:ruled}

In this appendix we review, and sketch the proofs, of several tools
from algebraic geometry that are required in our analysis, the main one of which is Theorem~\ref{th:ruledgz}.
These tools are presented in Guth and Zahl~\cite[Section 6]{GZ}, but
we reproduce them here, in a somewhat sketchy form, for the convenience
of the reader and in the interest of completeness. We work over the complex
field $\cplx$, but the results here also apply to our setting over the real
numbers (see \cite{GZ,SS4d} and a preceding remark for discussions of this issue).

A subset of $\cplx^N$ described by some polynomial equalities and one
non-equality, of the form
$$
\{p \in \cplx^N \mid f_1(p)=0, \ldots, f_r(p)=0, g(p)\ne 0\},
\quad \text{for\;} f_1,\ldots,f_r, g \in \cplx[x_1,\ldots, x_N] ,
$$
is called \emph{locally closed}. We recall that the (geometric) degree of an algebraic
variety $V\subset \cplx^N$ is defined as the number of intersection points of $V$ with
the intersection of $N-\dim(V)$ hyperplanes in general position
(see, e.g., Harris~\cite[Definition 18.1]{Harris}).
Locally closed sets have the following property.
\begin {theorem}[B\'ezout's inequality; B\"urgisser et al.~\protect{\cite[Theorem 8.28]{BCS}}] \label{th:bein}
Let $V$ be a nonempty locally closed set in $\cplx^N$, and let $H_1,\ldots, H_r$ be algebraic
hypersurfaces in $\cplx^N$. Then
$$
\deg(V\cap H_1\cap \cdots \cap H_r) \le \deg(V)\cdot \deg(H_1)\cdots \deg(H_r) .
$$
\end {theorem}

A constructible set $\C$ is easily seen to be a union of locally closed sets.
Moreover, one can decompose $\C$ uniquely as the union of irreducible locally closed sets
(namely, sets that cannot be written as the union of two nonempty and distinct locally closed sets).
By~B\"urgisser et al.~\cite[Definition 8.23]{BCS}), the degree of $\C$ is the
sum of the degrees of its irreducible locally closed components. Theorem~\ref{th:bein}
implies that when a constructible set
$\C$ has complexity $O(1)$, its degree is also $O(1)$. We
also have the following corollaries.
\begin {corollary} \label{co:bez}
Let
$$
X = \{p \in \cplx^N \mid f_1(p)=0, \ldots, f_r(p)=0, g(p)\ne 0\},
\quad \text{for\;} f_1,\ldots,f_r, g \in \cplx[x_1,\ldots, x_N] .
$$
If $X$ contains more than $\deg(f_1)\cdots \deg(f_r)$ points, then $X$ is infinite.
\end {corollary}
\noindent{\bf Proof.}
Assume that $X$ is finite. Put $V=\{p\in \cplx^n \mid g(p)\ne 0\}.$
By B\'ezout's inequality (Theorem~\ref{th:bein}), we have
$$\deg(X)\le \deg(V) \cdot \deg(f_1)\cdots \deg(f_r) = \deg(f_1)\cdots \deg(f_r),$$
where the equality $\deg(V)=1$ follows by the definition of the
degree of locally closed sets (see, e.g., B\"urgisser et
al.~\cite[Definition 8.23]{BCS}). When $X$ is finite, i.e.,
zero-dimensional, its degree is equal to the number of points in it,
counted with multiplicities. This implies that the number of points
in $X$ is at most $\deg(f_1)\cdots \deg(f_r)$, contradicting the
assumption of the theorem. Therefore, $X$ is infinite. \proofend

As an immediate consequence, we also have:
\begin {corollary} \label{co:const}
Let $\C\subset \cplx^N$ be a constructible set and write it as the union
of locally closed sets $\bigcup_{i=1}^t X_i$, where
$$
X_i = \{p \in \cplx^N \mid f^i_1(p)=0, \ldots, f^i_{r_i}(p)=0, g^i(p)\ne 0\},
\quad \text{for\;} f^i_1,\ldots,f^i_{r_i}, g^i \in \cplx[x_1,\ldots, x_N] .
$$
If $\C$ contains more than $\sum_{i=1}^t \deg(f^i_1)\cdots \deg(f^i_{r_i})$ points, then $\C$ is infinite.
\end {corollary}
For a constructible set $\C$, let $d(\C)$ denote the minimum of
$\sum_{i=1}^t \deg(f^i_1)\cdots \deg(f^i_{r_i})$, as in Corollary~\ref{co:const},
over all possible decompositions of $\C$ as the union of locally closed sets.
By B\'ezout's inequality (Theorem~\ref{th:bein}), it follows that $\deg(\C)\le d(\C)$.
Corollary~\ref{co:const} implies that if $\C$ contains more than $d(\C)$ points, then it is infinite.

Following Guth and Zahl~\cite[Section 4]{GZ}, we call an algebraic
curve $\gamma\subset \cplx^3$ a \emph{complete intersection} if
$\gamma=Z(P,Q)$ for some pair of polynomials $P,Q$. We let
$\cplx[x,y,z]_{\le E}$ denote the space of complex trivariate
polynomials of degree at most $E$, and choose an identification of
$\cplx[x,y,z]_{\le E}$ with\footnote{%
  Here $\binom {E+3} 3$ is the maximum number of monomials of the
  polynomials that we consider. For obvious reasons, the actual
  representation should be in the complex projective space
  $\cplx \P^{\binom {E+3} 3}$, but we use the many-to-one representation in
  $\cplx^{\binom {E+3} 3}$ for convenience.}
$\cplx^{\binom {E+3} 3}$.  We use the variable $\alpha$ to denote an
element of $(\cplx[x,y,z]_{\le E})^2$, and write
$$
\alpha = (P_\alpha, Q_\alpha)\in \left(\cplx[x,y,z]_{\le E}\right)^2
= \left(\cplx^{\binom {E+3} 3}\right)^2.
$$
Given an irreducible curve $\gamma$, we associate with it a choice
of $\alpha \in \left(\cplx^{\binom {E+3} 3}\right)^2$ such that
$\gamma$ is contained in $Z(P_\alpha,Q_\alpha)$, and the latter is a
curve (one can show that such an $\alpha$ always exists; see
Guth and Zahl~\cite[Lemma 4.2]{GZ} and also Basu and Sombra~\cite{BS}).
Let $x\in \gamma$ be a non-singular point\footnote{%
  Given an irreducible curve in $\reals^3$, a point $x\in \gamma$ is non-singular
  if there are polynomials $f_1, f_2$ that vanish on $\gamma$ such that
  $\nabla f_1(x)$ and $\nabla f_2(x)$ are linearly independent.}
of $\gamma$; we say that $\alpha$ is
associated to $\gamma$ at $x$, if $\alpha$ is associated to
$\gamma$, and $\nabla P_\alpha(x)$ and $\nabla Q_\alpha(x)$ are
linearly independent. We refer the reader to \cite[Definition 4.1 and Lemma 4.2]{GZ}
for details. This is analogous to the works of Guth and Katz~\cite{GK2} and
of Sharir and Solomon~\cite{SS4d} for the special cases of parameterizing lines in
three and four dimensions, respectively.

Before applying this machinery to derive the main result of this appendix,
we fill in the missing part in the proof of Theorem~\ref{incth}.
\begin{lemma} \label{c-const}
Let $V$ be some irreducible algebraic surface in $\cplx^3$, and let $\F$
be a constructible family of algebraic surfaces in $\cplx^3$ of constant
degree $E$. Let $\C$ be the family of the irreducible components of the
intersection curves of pairs of surfaces in $\F$.
Then $\C$ is a constructible family of curves.
\end{lemma}
\noindent{\bf Proof.}
We use the preceding identification of $\cplx[x,y,z]_{\le E}$ with
$\cplx^{\binom {E+3} 3}$, and use the variable $\alpha$ to denote an
element of $(\cplx[x,y,z]_{\le E})^2$, writing, as before
$$
\alpha = (P_\alpha, Q_\alpha)\in \left(\cplx[x,y,z]_{\le E}\right)^2
= \left(\cplx^{\binom {E+3} 3}\right)^2.
$$
Since $\F$ is constructible, the following set
$$
W = \{ \alpha = (P_\alpha, Q_\alpha)\in \left(\cplx[x,y,z]_{\le E}\right)^2
\mid P_\alpha, Q_\alpha \in \F\}
$$
is also contructible. By definition, recalling that $\C_{3,D}$ denotes the set of
all irreducible algebraic curves of degree at most $D$ in $\cplx^2$, we have
$$
\C = \{ \gamma \in \C_{3,E^2} \mid \exists \alpha \in W , \gamma \in P_\alpha \cap Q_\alpha\}.
$$
It then follows that $\C$ is irreducible by a simple variant of Lemma 4.3~of
Guth and Zahl~\cite{GZ}.
$\Box$

We now go on to the main result of the appendix.
In what follows, we fix a constructible set $\C_0 \subset \C_{3,E}$ of
irreducible curves of degree at most $E$ in 3-dimensional space (recall
that the entire family $\C_{3,E}$ is constructible). Following
\cite[Section 9]{GZ}, we call a point $p\in Z(f)$, for a given polynomial
$f\in\cplx[x,y,z]$, a \emph{$(t, \C_0, r)$-flecnode}, if there are
at least $t$ curves $\gamma_1,\ldots, \gamma_t \in \C_0$, such that, for
each $i=1,\ldots, t$, (i) $\gamma_i$ is incident to $p$, (ii) $p$ is a
non-singular point of $\gamma_i$, and (iii) $\gamma_i$ osculates to $Z(f)$
to order $r$ at $p$. This is a generalization of the notion of a
\emph{flecnodal} point, due to Salmon~\cite[Chapter XVII, Section III]{salmon}
(see also ~\cite{GK2, SS4d} for details).

With all this machinery, we can now present a (sketchy) proof of
Theorem~\ref{th:ruledgz}. The theorem is stated in Section~\ref{sec:intro},
and we recall it here. It is adapted from Guth and Zahl~\cite[Corollary 10.2]{GZ},
serves as a generalization of the Cayley--Salmon theorem on surfaces ruled by lines
(see, e.g., Guth and Katz~\cite{GK2}), and is closely related to Theorem~\ref{th:gz}
(also due to Guth and Zahl~\cite{GZ}).

\smallskip

\noindent{\bf Theorem \XX.}
{\it For given integer parameters $c,E$, there are constants $c_1=c_1(c,E)$,
$r=r(c,E)$, and $t=t(c,E)$, such that the following holds. Let $f$ be a complex
irreducible polynomial of degree $D \gg E$, and let $\C_0 \subset \C_{3,E}$ be
a constructible set of complexity at most $c$. If there exist at least $c_1 D^2$
curves of $\C_0$, such that each of them is contained in $Z(f)$ and contains at
least $c_1 D$ points on $Z(f)$ that are $(t,\C_0,r)$-flecnodes, then $Z(f)$ is
infinitely ruled by curves from $\C_0$. In particular, if $Z(f)$ is not infinitely
ruled by curves from $\C_0$ then, except for at most $c_1 D^2$ \emph{exceptional}
curves, every curve in $\C_0$ that is fully contained in $Z(f)$ is incident to at
most $c_1 D$ points that are incident to at least $t$ curves in $\C_0$ that are also
fully contained in $Z(f)$.}

\smallskip

\noindent {\bf Proof.}
For the time being, let $r$ be arbitrary.
By Guth and Zahl~\cite[Lemma 8.3 and Equation (8.1)]{GZ}, since $f$ is irreducible,
there exist $r$ polynomials\footnote{%
  To say that $h_j$ is a ploynomial in $\alpha$ (and $p$) means that it is a polynomial in the
  $2\binom{E+3}{3}$ coefficients of the monomials of the two polynomials in the pair $\alpha$
  (and in the coordinates $(x,y,z)$ of $p$).}
$h_j(\alpha,p)\in \cplx[\alpha,x,y,z]$, for $j=1,\ldots, r$ of
degree at most $b_j$ in $\alpha$ (where $b_j$ is a constant
depending on $j$ and on $E$), and of degree $O(D)$ in $p=(x,y,z)$,
with the following property: let $\gamma$ be an irreducible curve,
let $p$ be a non-singular point of $\gamma$, and let $\alpha$ be
associated to $\gamma$ at $p$, then $\gamma$ osculates to $Z(f)$ to
order $r$ at $p$ if and only if $h_j(\alpha,p)=0$ for $j=1,\ldots,
r$. (These polynomials are suitable representations of the first $r$
terms of the Taylor expansion of $f$ at $p$ along $\gamma$;
see~\cite[Section 6.2]{GZ} for this definition, and
also~\cite{GK2,SS4d} for the special cases of lines in $\reals^3$
and $\reals^4$, respectively.)

Regarding $p$ as fixed, the system $h_j(\alpha,p)=0$, for $j=1,\ldots,r$, in
conjunction with the constructible condition that $\alpha \in \C_0$, defines
a constructible set $\C_p$. By definition, we have
$d(\C_p) \le \left(\prod_{j=1}^r b_j\right) \cdot d(\C_0)$, which is a constant that
depends only on $r$ and $E$. By Corollary~\ref{co:const}, $\C_p$ is either infinite or
contains at most $d(\C_p)=O(1)$ points. By Guth and Zahl~\cite[Corollary 12.1]{GZ},
there exist a Zariski open set $\O$, and a sufficiently large constant $r_0$, that depend
on $\C_0$ and $E$ (see~\cite[Theorem 8.1]{GZ} for the way $r_0$ is obtained),
such that if $p\in\O$ is a $(t,\C_0, r)$-flecnode, with $r\ge r_0$,
there are at least $t$ curves that are incident to $p$ and are fully
contained in $Z(f)$. Since, by assumption, there are at least $c_1
D^2$ curves, each containing at least $c_1 D$ $(t, \C_0,
r)$-flecnodes, it follows from~\cite[Proposition 10.2]{GZ} that
there exists a Zariski open subset ${\O}$ of $Z(f)$, all of whose
points are $(t,\C_0,r)$-flecnodes. As noted above,~\cite[Corollary
12.1]{GZ} then implies that every point of ${\O}$ is incident to at
least $t$ curves of degree at most $E$ that are fully contained in
$Z(f)$. As observed above, when $t\ge \left(\prod_{j=1}^r b_j\right)
\cdot d(\C_0)$, a constant depending only on $\C_0$ and $E$, $Z(f)$
is infinitely ruled on this Zariski open set. By a simple argument
(a variant of which is given in~\cite[Lemma 6.1]{SS4dv}), we can
conclude that $Z(f)$ is infinitely ruled (everywhere) by curves from
$\C_0$, thus completing the proof. \proofend


\end{document}